\documentclass[11pt,a4paper]{article}

\usepackage{amsmath,amssymb}
\newcommand{\e}{\varepsilon}
\newcommand{\D}{\Delta}
\newcommand{\La}{\Lambda}
\newcommand{\va}{\varphi}
\newcommand{\n}{\nabla}
\newcommand{\N}{\frac{N}{2}}

\newcommand{\g}{\int_{\mathbb{R}^{N}}}
\newcommand{\p}{\partial}

\newcommand{\R}{\mathbb{R}}
\newcommand{\h}{\hookrightarrow}

\newtheorem{definition}{Definition}
\newtheorem{theorem}{Theorem}
\newtheorem{notation}{Notation}
\newtheorem{proposition}{Proposition}
\newtheorem{corollary1}{Corollary}
\newtheorem{remarka}{Remark}

\newtheorem{lemme}{Lemma}

\title{Existence of weak solution for compressible fluid models of Korteweg type }
\author{Boris Haspot \thanks{Karls Ruprecht Universit\"at Heidelberg, Institut for Applied Mathematics,
Im Neuenheimer Feld 294,
D-69120 Heildelberg, Germany.
Tel. 49(0)6221-54-6112
}}
\date{}
\begin{document}
\maketitle
\begin{abstract}
This work is devoted to proving existence of global weak solutions
for a general isothermal model of capillary fluids derived by J.- E
Dunn and
J. Serrin (1985) \cite{3DS}, which can be used as a phase transition model.\\
We  improve the results of \cite{3DD} by showing the existence of global weak solution in dimension two for initial data in the energy space, close to a stable equilibrium and with specific choices
on the capillary coefficients. In particular we are interested in capillary coefficients approximating a constant capillarity coefficient $\kappa$.
To finish we show the existence of global weak solution in dimension one for a specific type of capillary coefficients
with large initial data in the energy space.
\end{abstract}
\section{Introduction}
\subsection{Derivation of Korteweg model}
We are concerned with compressible fluids endowed with internal
capillarity. The model we consider  originates from the XIXth
century work by van der Waals and Korteweg \cite{3K} and was
actually derived in its modern form in the 1980s using the second
gradient theory, see for instance \cite{3JL,3TN}.
\\
Korteweg-type models are based on an extended version of
nonequilibrium thermodynamics, which assumes that the energy of the
fluid not only depends on standard variables but also on the
gradient of
the density.
Let us now consider a fluid of density $\rho\geq 0$, velocity field $u\in\R^{N}$, we are now interested in the following
compressible capillary fluid model, which can be derived from a Cahn-Hilliard like free energy (see the
pioneering work by J.- E. Dunn and J. Serrin in \cite{3DS} and also in
\cite{3A,3C,3GP}).
The conservation of mass and of momentum write:
\begin{equation}
\begin{cases}
\begin{aligned}
&\frac{\p}{\p t}\rho+{\rm div}(\rho u)=0,\\
&\frac{\p}{\p t}(\rho u)+{\rm div}(\rho
u\otimes u)-\rm div(2\mu(\rho) D (u))-\n\big(\lambda(\rho)){\rm div}u\big)+\n P(\rho)={\rm div}K,
\end{aligned}
\end{cases}
\label{3systeme}
\end{equation}
where the Korteweg tensor read as following:
\begin{equation}
{\rm div}K
=\n\big(\rho\kappa(\rho)\D\rho+\frac{1}{2}(\kappa(\rho)+\rho\kappa^{'}(\rho))|\n\rho|^{2}\big)
-{\rm div}\big(\kappa(\rho)\n\rho\otimes\n\rho\big).
\label{divK}
\end{equation}
$\kappa$ is the coefficient of capillarity and is a regular function. The term
${\rm div}K$  allows to describe the variation of density at the interfaces between two phases, generally a mixture liquid-vapor. $P$ is a general increasing pressure term.
$D (u)=\frac{1}{2}[\n u+^{t}\n u]$ being the stress tensor, $\mu$ and $\lambda$ are the two Lam\'e viscosity coefficients depending on the density $\rho$) and satisfying:
$$\mu>0\;\;\mbox{and}\;\;2\mu+N\lambda\geq0.$$
Here we want to investigate the existence of global weak solution for the system (\ref{3systeme}), so we have to describe precisely the form of the non linear terms coming from the capillary tensor.
In this goal one can now rewrite $K$.For this and 
to simplify the presentation, we assume only that $\kappa(\rho)=\kappa\rho^{\alpha}$. With this choice we have to distinguish the case $\alpha\ne -2$ and $\alpha=-2$.
We refer to the appendix in section \ref{appendix} for more details on the formal computations on $K$.
We get then:
\begin{equation}
\begin{aligned}
&K_{i,j}=\big(A^{1}_{\alpha}\D B(\rho)-A^{2}_{\alpha}
|\n A(\rho)|^{2}\,
\big)\p_{i,j}-B_{\alpha}\p_{i}A(\rho)\p_{j}A(\rho)\;\;\;\;&\mbox{if}\;\;\;\alpha\ne -2,\\
&K_{i,j}=\kappa\big(\D \log(\rho)+\frac{1}{2}|\n \log(\rho)|^{2}\,
\big)\p_{i,j}-\kappa\p_{i} \log(\rho)\p_{j} \log(\rho)\;\;\;\;\;\;&\mbox{if}\;\;\;\alpha=-2.\\
\end{aligned}
\label{divK1}
\end{equation}
with:
$$
\begin{aligned}
&A(\rho)=\rho^{\frac{\alpha}{2}+1},\;\;\;B(\rho)=\rho^{2+\alpha}\;\;\;\mbox{and}\;\;\;\;A^{1}_{\alpha}=
\frac{\kappa}{2+\alpha},\;\;\;\;
A^{2}_{\alpha}=\frac{2\kappa(\alpha+1)
}{(\alpha+2)^{2}},\;\;\;\;B_{\alpha}=\frac{4\kappa}{(\alpha+2)^{2}}.
\end{aligned}
$$
For the general case, we have $A^{'}(\rho)=\sqrt{\kappa}$ and $B^{'}(\rho)=\rho\kappa(\rho)$ (see the appendix for more details).
We remark then that the form of the non linear terms appearing in the tensor $K$ corresponds to quadratic gradient terms
$\p_{i}A(\rho)\p_{j}A(\rho)$ and some terms in $B(\rho)$ of pressure type. It is the main difference between the system of Korteweg and this one of Navier-Stokes compressible where the difficulty state in the treatment of the pressure term (see \cite{3F}, \cite{3L2}))
Before getting into the heart of mathematical results, one can recall
first derive the physical energy bounds of the Korteweg
system. Let $\bar{\rho}>0$ be a constant reference density, and let
$\Pi $ be defined
by:
$$\Pi(s)=s\biggl(\int^{s}_{\bar{\rho}}\frac{P(z)}{z^{2}}dz-\frac{P(\bar{\rho})}{\bar{\rho}}\biggl),$$
so that $P(s)=s\Pi^{'}(s)-\Pi(s)\, ,\,\Pi^{'}(\bar{\rho})=0$
and:
$$\p_{t}\Pi(\rho)+{\rm div}(u\Pi(\rho))+P(\rho){\rm div}(u)=0\;\;\;\mbox{in}\;\;
{\cal D}^{'}((0,T)\times\R^{N}).$$
Notice that $\Pi$ is convex as far as $P$ is non decreasing (since
$P^{'}(s)=s\Pi^{''}(s)$), which is the case for $\gamma$-type
pressure laws.
Multiplying the equation of momentum conservation in the system
(\ref{3systeme}) by $u$ and integrating by parts over $\R^{N}$,
we obtain the following
estimate ( for more details on the integration by parts, especially on the term ${\rm div}K$ we refer to the appendix
section \ref{appendix}):
\begin{equation}
\begin{aligned}
&\int_{\R^{N}}\big(\frac{1}{2}\rho
|u|^{2}+(\Pi(\rho)-\Pi(\bar{\rho}))+\frac{1}{2}\kappa(\rho)|\nabla\rho|^{2}\big)(t)dx
+\int_{0}^{t}\int_{\R^{N}}\big(\mu(\rho)
|D(u)|^{2}\\
&\hspace{0,5cm}+\xi(\rho)|{\rm div} u|^{2}\big)dx
\leq\int_{\R^{N}}\big(\frac{|m_{0}|^{2}}{2\rho}+(\Pi(\rho_{0})-\Pi(\bar{\rho}))
+\frac{\kappa(\rho_{0})}{2}|\nabla\rho_{0}|^{2}\big)dx.
\label{3inegaliteenergie1}
\end{aligned}
\end{equation}
with $\xi(\rho)=\mu(\rho)+\lambda(\rho)$.
We will note in the sequel:
\begin{equation}
{\cal E}(t)=\int_{\R^{N}}\big(\frac{1}{2}\rho
|u|^{2}+(\Pi(\rho)-\Pi(\bar{\rho}))+\frac{\kappa(\rho)}{2}|\nabla\rho)|^{2}\big)(t)dx,
\label{3defenergie}
\end{equation}
It follows that assuming that the initial total energy is finite:\\
$${\cal E}_{0}=\int_{\R^{N}}\big(\frac{|m_{0}|^{2}}{2\rho}+(\Pi(\rho_{0})-\Pi(\bar{\rho}))
+\frac{\kappa(\rho_{0})}{2}|\nabla\rho_{0}|^{2}\big)dx<+\infty\,,$$ then we
have the a priori following bounds:
$$\Pi(\rho)-\Pi(\bar{\rho}),\;\;\mbox{and}\;\;\rho |u|^{2}\in L^{1}(0,\infty,L^{1}(\R^{N})),$$
$$\sqrt{\kappa(\rho)}\n\rho\in L^{\infty}(0,\infty,L^{2}(\R^{N}))^{N},\;\;\mbox{and}\;\;\n u\in
L^{2}(0,\infty,\R^{N})^{N^{2}}.$$
In the sequel, we aim at solving
the problem of global existence of weak solution for the system
(\ref{3systeme})so  
assuming that we dispose from smooth approximates sequel solutions $(\rho_{n},u_{n})_{n\in\mathbb{N}}$ of system (\ref{3systeme}).
One can remark then easily that the main difficulty lies in the fact to be able to pass to the limit in the quadratic term
$\n A(\rho_{n})\otimes\n A(\rho_{n})$ which belongs only to $L^{\infty}(L^{1})$.
According to the classical theorems on weak topology,
$\n A(\rho_{n})\otimes\n A(\rho_{n})$ converges up to extraction to a
measure $\nu$, so how we can prove that
$\nu=\n A(\rho)\otimes\n A(\rho)$ where $\rho$ is the limit of the sequence
$(\rho_{n})_{n\in\mathbb{N}}$
in appropriate space.
Notice that if we compare the Korteweg system with compressible Navier-Stokes system, the main difficulty for proving global existence of weak solution ( see \cite{3L2}, \cite{3F}) consists to find strong compactness properties on the density $\rho$ in $L^{p}_{loc}$ spaces  to pass to the limit in the pressure term $P(\rho)=a\rho^{\gamma}$ with $\gamma>\frac{N}{2}$, $N\geq 2$ in the case of the works of E. Feireisl. In the capillary case, more a priori bounds are available for the density, as $\n A(\rho)$
belongs to $L^{\infty}(H^{1}(\R^{N}))$. Hence in our case in dimension $N=2$, one can easily pass to the limit in the pressure term.\\
Another difficulty in compressible fluid mechanics is to
deal with the vacuum and we will see that this problem does appear
in the model of Korteweg, when we want get some estimates on $\n\rho$. As a
matter of fact, the existence of global weak solution
for the model of Korteweg with constant capillary coefficient for $N\geq 2$
is still an open problem even in the case $N=2$ if we don't assume some condition on the vacuum as we will see in the sequel.
The first  ones to have studied the problem, are R. Danchin and B.
Desjardins in \cite{3DD}. They showed that if we take initial data
close to a stable equilibrium in the energy space and assume that we
control the vacuum  and the norm $L^{\infty}$ of the density $\rho$,
then we get global weak solution
 in dimension two. Controlling the vacuum amounts here to get a bound of $\frac{1}{\rho}$ in $L^{\infty}$.
Recently D. Bresch, B. Desjardins and C-K. Lin in \cite{3BDL} got
some global weak solutions for the isotherm Korteweg model with some
specific viscosity coefficients. Indeed, they assume that
$\mu(\rho)=C\rho$ with $C>0$ and $\lambda(\rho)=0$. By choosing
these specific coefficients they obtain a gain of derivatives on the
density $\rho$ where $\rho$ belongs to $L^{2}(H^{2})$.
It is easy at that time with this type of estimate on the density $\rho$ to get strong convergence on the term of capillarity. However a new difficulty takes place
concerning the loss of information on the gradient of $u$ when vacuum appearing (indeed the viscosity coefficients are degenerated), so it becomes
involved to pass to the limit in the term $\rho_{n} u_{n}\otimes u_{n}$. That's why
the solutions of D. Bresch, B. Desjardins and C-K. Lin require some
specific test
functions which depend on the density $\rho$. This test functions depending of the solution allow to deal with the vacuum.
\\
Let us mention briefly that the existence of strong solutions for $N\geq2$ is known since the works by H. Hattori and D. Li \cite{3H1,3H2}.
Notice that high order regularity in Sobolev space $H^{s}$ is required, namely the initial data $(\rho_{0},\rho_{0}u_{0})$
are assumed to belong to $H^{s}\times H^{s-1}$ with $s\geq\N+4$. Moreover they consider convex pressure profiles, which cannot cover the case of Van der Waals' equation of state. R. Danchin and B. Desjardins in \cite{3DD} improve this result by working in critical spaces for the scaling of the equations, more precisely the initial data $(\rho_{0},\rho_{0}u_{0})$
belong to $B^{\N}_{2,1}\times B^{\N-1}_{2,1}$ (the fact that $B^{\N}_{2,1}$ is embedded in $L^{\infty}$ play a crucial role to control the vacuum but or to work in multiplier space). In \cite{3MK}, M. Kotschote show the existence of strong solution for the isothermal model in bounded domain by using Dore\^a-Venni Theory and $\mathcal{H}^{\infty}$ calculus.
In \cite{3H}, we generalize the results of \cite{3DD} in the case of non isothermal Korteweg system with physical coefficients depending on the density and the temperature. We get strong solutions with initial data belonging to the critical spaces $B^{\N}_{2,1}\times B^{\N-1}_{2,1}\times B^{\N-2}_{2,1}$ when the physical coefficients depend only on the density.
When they depend on the temperature too, it requires more regular initial data to control the temperature. \\
Our present result improve the results of R. Danchin and
B. Desjardins in \cite{3DD} and D. Bresch, B. Desjardins an C-K. Lin in \cite{3BDL}, indeed we show the existence of global weak solution with small initial data in the energy space
for specific choices of the capillary coefficients and with general viscosity coefficient.
Comparing with the results of \cite{3BDL}, we get global weak solutions with test function $\va\in C^{0}_{\infty}(\R^{N})$ and not depending of the density $\rho$. Moreover our result can be applied for general viscosity coefficients.
Moreover in the case of a constant capillary coefficient, we show that we don't need
to control $\rho$
in $L^{\infty}$ norm to get global weak solution as in \cite{3DD}.
In fact we have extracted of the structure of capillarity term a new energy inequality using fractionary derivative which allows a
gain of derivative on the density $\rho$. This method  shows that the case $\kappa(\rho)=\frac{\kappa}{\rho^{2}}$, $\kappa>0$
is a critical case where we can get global weak solutions with small initial data in the energy spaces but without any conditions on the density as in \cite{3DD}. Moreover we generalize the existence of global weak strong solutions to the case of capillarity coefficients approximating the constant case. To finish we explain why it seems a little bit trick to get global weak solution with large initial data in the energy space, indeed the problem is then to control some possible concentration effect on $|\n A(\rho)|^{2}$.\\
In section \ref{3S1} we recall some definitions on the Orlicz space
and some classical energy inequalities on the system in these spaces. In section
\ref{3S2} we show a theorem of global existence of weak
solutions in dimension two for specific choices of capillarity coefficient with general viscosity coefficients. More precisely
we investigate a new structure on the capillarity coefficient which allow us to get a gain of derivative on $A(\rho)$ under
the condition to get for $\e>0$ enough small $\|1_{B(x,r)}\n A(\rho)\|_{L^{\infty}(L^{2})}\leq \e$ uniformly for all $x\in\R$ with
$r>0$ enough small. By imposing small conditions on the initial data, we obtain our results and we explain what happen with large initial data.
In the last section we investigate the case
of the dimension one, and we get a theorem of  existence of
global weak solution in the energy space with
large initial data.
\section{Classical a priori estimates and Orlicz spaces}
\label{3S1}
\subsection{Classical a priori estimates}
We first want recall a
priori bounds for initial density in Orlicz space, it means that we will work around a constant state $\bar{\rho}>0$.
We rewrite the mass equation by using renormalized solutions,
and the momentum equation.In the sequel to simplify we let $P(\rho)=a\rho^{\gamma}$ with $\gamma\geq 1$ and $a$ a positive regular function. We get the following formal identities:
\begin{equation}
\begin{cases}
\begin{aligned}
&\frac{1}{\gamma-1}\frac{\p}{\p
t}\big(\rho^{\gamma}-\bar{\rho}^{\gamma}
-\gamma\bar{\rho}^{\gamma-1}(\rho-\bar{\rho})\big)+{\rm
div}\big[u\frac{\gamma}{\gamma-1}(\rho^{\gamma}-\bar{\rho}^{\gamma-1}\rho)\big]=u\cdot\n(\rho^{\gamma})\\
&\rho\frac{\p}{\p t}\frac{|u|^{2}}{2}+\rho
u\cdot\n\frac{|u|^{2}}{2}-\rm div(2\mu(\rho) D (u))\cdot u-\n\big(\lambda(\rho)){\rm div}u\big)\cdot u+au\cdot\n\rho^{\gamma}\\
&\hspace{10,6cm}=u\cdot{\rm div}K\;,
\label{3a26}
\end{aligned}
\end{cases}
\end{equation}
\begin{notation}
In the sequel we will note:
$$j_{\gamma}(\rho)=\rho^{\gamma}+(\gamma-1)\bar{\rho}^{\gamma}
-\gamma\bar{\rho}^{\gamma-1}\rho.$$
\end{notation}
Therefore we find by summing the two equalities of (\ref{3a26}):
\begin{equation}
\begin{aligned}
&\frac{\p}{\p
t}\big[\rho\frac{|u|^{2}}{2}+\frac{a}{\gamma-1}j_{\gamma}(\rho)\big]+{\rm
div}\big[u(\frac{a\gamma}{\gamma-1}(\rho^{\gamma}-\bar{\rho}^{\gamma-1}\rho)+\rho\frac{|u|^{2}}{2})\big]\\
&\hspace{4cm}-\rm div(2\mu(\rho) D (u))\cdot u
-\n\big(\lambda(\rho)){\rm div}u\big)\cdot u=u\cdot{\rm div}K.
\label{3a27}
\end{aligned}
\end{equation}
We may then integrate in space the equality (\ref{3a27}) and we get (see for more details the appendix \ref{appendix}):
\begin{equation}
\begin{aligned}
&\int_{\R^{N}}\big(\rho\frac{|u|^{2}}{2}+\frac{a}{\gamma-1}j_{\gamma}(\rho)+\frac{1}{2}\kappa(\rho)|\n\rho|^{2}\big)(t,x)dx
+\int_{0}^{t}ds\int_{\R^{N}}2\mu(\rho)|D
u|^{2}dx\\
&+\int_{0}^{t}ds\int_{\R^{N}}\lambda(\rho)|{\rm
div}u|^{2}dx\leq\int_{\R^{N}}\big(\rho_{0}\frac{|u_{0}|^{2}}{2}+
\frac{a}{\gamma-1}j_{\gamma}(\rho_{0})+\kappa(\rho_{0})|\n\rho_{0}|^{2}\big)(x)\,dx.
\end{aligned}
\label{3a29}
\end{equation}
\begin{notation}
In the sequel we will note:
$$
\begin{aligned}
&{\cal E}^{\gamma}(t)=\int_{\R^{N}}\big(\rho\frac{|u|^{2}}{2}+\frac{a}{\gamma-1}j_{\gamma}(\rho)+\frac{1}{2}\kappa(\rho)|\n\rho|^{2}\big)(t,x)dx,\\
&{\cal E}^{\gamma}_{0}=
\int_{\R^{N}}\big(\rho_{0}
\frac{|u_{0}|^{2}}{2}+\frac{a}{\gamma-1}j_{\gamma}(\rho_{0})+\frac{1}{2}\kappa(\rho_{0})|\n\rho_{0}|^{2}\big)dx.
\end{aligned}
$$
\end{notation}
We now want to estimate this quantity $j_{\gamma}(\rho)$ and in this
goal we recall some properties of Orlicz spaces.
\subsection{Orlicz spaces}
We begin by describing the Orlicz space in which we will work:
$$L^{q}_{p}(\R^{N})=\{f\in L^{1}_{loc}(\R^{N})/f 1_{\{|f|\leq\delta\}}\in L^{p}(\R^{N}),\;\;
f 1_{\{|f|\geq\delta\}}\in L^{q}(\R^{N})\}$$
where $\delta$ is fixed, $\delta>0$.\\
First of all, it is not difficult to check that $L^{q}_{p}$ does not
depend on the choice of $\delta>0$ since $\frac{x^{p}}{x^{q}}$ is
bounded from above and from below on any interval
$[\delta_{1},\delta_{2}]$ with $0<\delta_{1}\leq\delta_{2}<+\infty$.
In particular we deduce that we have:
$$f^{\e}\in L^{\frac{q}{\e}}_{\frac{p}{\e}}(\R^{N})\;\;\;\mbox{if}\;\;f\in L^{q}_{p}(\R^{N})\;\;\mbox{and}
\;\;p,\,q\geq\e .$$ Obviously we get
$\mbox{meas}\{|f|\geq\delta\}<+\infty$ if $f\in L^{q}_{p}(\R^{N})$
and thus we have the embedding:
$$L^{q}_{p}(\R^{N})\subset L^{q_{1}}_{p_{1}}(\R^{N})\;\;\;\mbox{if}\;\;1\leq q_{1}\leq q<+\infty,\;\;
1\leq p\leq p_{1}<+\infty. $$ Next, we choose $\Psi$ a convex
function on $[0,+\infty)$ which is equal (or equivalent) to $x^{p}$
for $x$ small and to $x^{q}$ for $x$ large, then we can define the
space $L^{q}_{p}(\R^{N})$ as follows:
\begin{definition}We define then the Orlicz space
$L^{q}_{p}(\R^{N})$ as follows:\\\texttt{}
$$L^{q}_{p}(\R^{N})=\{f\in L^{1}_{loc}(\R^{N})/\Psi(f)\in
L^{1}(\R^{N})\}.$$
\end{definition}
We can check that $L^{q}_{p}(\R^{N})$ is a linear vector space. Now
we endow $L^{q}_{p}(\R^{N})$ with  a norm so that
$L^{q}_{p}(\R^{N})$ is a separable Banach space:
$$\|f\|_{L^{q}_{p}(\R^{N})}=\inf\{t>0/\;\;\Psi(\frac{f}{t})\leq 1\}.$$
We recall now some useful properties of Orlicz spaces.
\begin{proposition}
The following properties hold:
\begin{enumerate}
\item Dual space: If $p>1$ and $q>1$
then $(L^{q}_{p}(\R^{N}))^{'}=L^{q^{'}}_{p^{'}}(\R^{N})$ where
$q^{'}=\frac{q}{q-1},\,p^{'}=\frac{p}{p-1}$.
\item $L^{q}_{p}=L^{p}+L^{q}$ if $1\leq q\leq p<+\infty$.
\item Composition: Let $F$ be a continuous function on $\R$ such that $F(0)=0$, $F$
is differentiable at $0$ and $F(t)|t|^{-\theta}\rightarrow\alpha\ne
0$ at $t\rightarrow +\infty$. Then if $q\geq\theta$,
$$F(f)\in L^{\frac{q}{\theta}}_{p}(\R^{N})\;\;\mbox{if}\;\;f\in
L^{q}_{p}(\R^{N}).$$
\end{enumerate}
\end{proposition}
Now we can recall a property on the Orlicz space concerning the
inequality of energy.
\begin{proposition}
The function $j_{\gamma}(\rho)$ is in $ L^{1}(\R^{N})$ if and only
if $\rho-\bar{\rho}\in L^{\gamma}_{2}.$
\end{proposition}
{\bf Proof:}
On the set $\{|\rho-\bar{\rho}|\leq\delta\}$, $\rho$ is bounded from
above, since $\gamma>1$ we thus deduce that $j_{\gamma}(\rho)$ is
equivalent to $|\rho-\bar{\rho}|^{2}$ on the set
$\{|\rho-\bar{\rho}|\leq\delta\}$.
Next on the set  $\{|\rho-\bar{\rho}|\geq\delta\}$, we observe that
for some $\nu\in (0,1)$ and $C\in (1,+\infty)$, we have:
$$\nu|\rho-\bar{\rho}|^{\gamma}\leq j_{\gamma}(\rho)\leq C|\rho-\bar{\rho}|^{\gamma}.$$
\hfill {$\Box$}
\subsubsection*{Link with our energy estimate}
We recall the definition of the fractional derivative operator
$\La^{s}$.
\begin{definition}
We define the operator $\La^{s}$ as follows:
$\widehat{\La^{s} f}=|\xi|^{s}\widehat{f}$.
\end{definition}
We give now some useful results concerning Sobolev spaces, we start with a proposition
coming from a theorem of interpolation by Riesz-Thorin.
\begin{proposition}
The Fourier transform is continuous from $L^{p}$ in $L^{q}$ with
$p\in[1,2]$, $q\in[2,+\infty]$ and
$\frac{1}{p}+\frac{1}{q}=1.$
\end{proposition}
We recall here the definition of homogeneous Sobolev space.
\begin{definition}
Let $s\in\R$. $f$ is in the homogeneous space $\dot{H}^{s}$ if:
$|\xi|^{s}\hat{f}\in L^{2}(\R^{N})$.
\end{definition}
\begin{proposition}
Let $f\in\dot{H}^{s}$ with $s>0$ and $f\in L^{p}+L^{2}$ with $1\leq
p<2$. Then $f\in L^{2}$. \label{3P1}
\end{proposition}
{\bf Proof:}
Indeed we have as $f\in\dot{H}^{s}$:
$$\int_{\R^{N}}|\xi|^{2s}|\widehat{f}|^{2}d\xi<+\infty,$$
so $\widehat{f}1_{\{|\widehat{f}|\geq1\}}\in L^{2}(\R^{N})$. And as
$f=f_{1}+f_{2}$ with $f_{1}\in L^{p}(\R^{N})$ and $f_{2}\in L^{2}$.
By using the Riesz-Thorin theorem, we know that
$\widehat{f_{1}}\in L^{q}(\R^{N})$ with $\frac{1}{p}+\frac{1}{q}=1$.
As $q\geq2$ we then have $\widehat{f}1_{\{|\widehat{f}|\leq1\}}\in
L^{2}(\R^{N})$. This concludes the proof.
\hfill {$\Box$}\\
According to the above theorem and our energy estimate we get that
for all $T\in\R$,
$\rho-\bar{\rho}\in L^{\infty}(0,T;L^{\gamma}_{2}(\R^{N}))$.
\begin{remarka}
\label{3commentaire}
Let $\kappa(\rho)=\kappa\rho^{\alpha}$ with $\alpha\ne -2$
We assume now that $\frac{1}{\rho}\in L^{\infty}$ if $\alpha>-2$ and $\rho\in L^{\infty}$ if $\alpha<-2$
(These hypothesis will be supposed in the theorem \ref{3T1}).
We have then by using previous properties on Orlicz spaces and (\ref{3a29}):
\begin{itemize}
\item if $\gamma\geq2$ then $L^{\gamma}_{2}(\R^{N})\hookrightarrow L^{2}(\R^{N})$ and
so $\rho-\bar{\rho}\in L^{\infty}(H^{1}(\R^{N}))$.
\item if $\gamma\leq2$ then following the proposition \ref{3P1} and the fact that $L^{\gamma}_{2}=L^{\gamma}
+L^{2}$ we get $\rho-\bar{\rho}\in L^{\infty}(H^{1}(\R^{N}))$.
\end{itemize}
\end{remarka}
We finish with one proposition which give some results of refined Sobolev, for a proof see \cite{3BCD}.
\begin{proposition}
\label{3refined}
Let $1\leq q<p<+\infty$ and $\alpha$ be a positive real number. A constant $C$ exists such that:
$$\|f\|_{L^{p}}\leq C\|f\|^{1-\theta}_{B^{-\alpha}_{\infty,\infty}}\|f\|^{\theta}_{B^{\beta}_{q,q}}\;\;\;\mbox{with}
\;\;\beta=\alpha(\frac{p}{q}-1)\;\;\;
\mbox{and}\;\;\theta=\frac{q}{p}.$$
\end{proposition}
We end with a classical proposition of weak topology.
\begin{proposition}
Let $1<r<+\infty$ and $K$ a compact. Let $f_{k}$ a sequel of vector fields in $\R^{N}$ such that:
\begin{itemize}
\item $f_{k}$ is uniformly bounded in $L^{r}(K)$ and $f_{k}$ converge a.e to $f$,
\end{itemize}
then $f_{k}$ converge weakly to $f$ in $L^{r}(K)$.
\label{3propzuily}
\end{proposition}
{\bf Proof:} We have to show that for all $\va\in L^{r^{'}}$ with $\frac{1}{r^{'}}+\frac{1}{r}=1$ we have:
$$\big(f_{k},\va\big)\rightarrow_{k\rightarrow+\infty}\big(f,\va\big).$$
It is enough in fact to show this property for all $\va\in\xi$ with $\xi$ a dense space in $L^{r^{'}}$. By the Fatou theorem we check that $f\in L^{r}(K)$. Next we set for all $k\in\mathbb{N}$:
$$E(k)=\{y\in K:\;|f_{k}(y)-f(y)|>1\},$$
and $F_{p}=\bigcup_{k\geq p}E(k)$. We have now $N=\bigcap_{p\geq 0}F_{p}$ with:
$$N=\{y\in K:\;\;\forall p\;\exists k\geq p\;\;|f_{k}(y)-f(y)|>1\}$$
As $f_{k}$ converge a.e to $f$ we have the measure of Lebesgue of $N$ which is null.
We choose then $\xi=\{\va\in L^{r^{'}}(K):\;\;\mbox{supp}\va\subset K\backslash N\}$.
We show next easily by a convergence dominated that for all $\va\in\xi$, $(f_{k}-f,\va)\rightarrow_{k\rightarrow+\infty}0$. And this achieve
the proof.\hfill {$\Box$}\\
We can now explain what we mean by weak solution of problem
(\ref{3systeme}) in dimension $N=2$.
\begin{definition}
\label{3defexistence} Let the couple $(\rho_{0},u_{0})$ satisfy;
\begin{enumerate}
\item $\rho_{0}\in L^{\gamma}_{2}(\R^{N})$, $\n A(\rho_{0})\in L^{2}(\R^{N})$
\item $\rho_{0}|u_{0}|^{2}\in L^{1}(\R^{N})$
\item $\rho_{0}u_{0}=0$ whenever $x\in\{\rho_{0}=0\}$,
\end{enumerate}
We have the following definition:
\begin{itemize}
\item A couple $(\rho,u)$ is called a weak solution of problem
(\ref{3systeme}) on $I\times\R^{N}$ with $I$ an interval of $\R$ if:
\begin{itemize}
\item $\rho\in L^{\infty}(L^{\gamma}_{2}(\R^{N}))$, $\n A(\rho)\in
L^{\infty}(L^{2}(\R^{N}))$, $\va\rho\in L^{2}(H^{1+\alpha}(\R^{N}))$
$\forall\alpha\in]0,1[$, and $\forall\va\in C^{\infty}_{0}(\R^{N})$.
\item $\n u\in L^{2}(L^{2}(\R^{N}))$, $\rho|u|^{2}\in
L^{\infty}(L^{1})$.
\item Mass equation holds in ${\cal D}^{'}(I\times\R^{N})$.\label{3vrai1}
\item Momentum equation holds in ${\cal D}^{'}(I\times\R^{N})^{N}$.\label{3vrai2}
\item $\lim_{t\rightarrow
0^{+}}\int_{\R^{N}}\rho(t)\va=\int_{\R^{N}}\rho_{0}\va$,
$\forall\va\in{\cal D}(\R^{N})$,\label{3vrai3}
\item $\lim_{t\rightarrow
0^{+}}\int_{\R^{N}}\rho u(t)\cdot\phi=\int_{\R^{N}}(\rho
u)_{0}\cdot\phi$, $\forall\phi\in{\cal
D}(\R^{N})^{N}$.\label{3vrai4}
\end{itemize}
\item The quantity ${\cal E}^{\gamma}$ is finite and inequality (\ref{3a29}) holds a.e in $I$.
\end{itemize}
\end{definition}
\section{Existence of global weak solutions for $N=2$}
\label{3S2}
\subsection{Gain of derivatives in the case $N=2$}
We give now a general description of the viscosity coefficients and
in the sequel we will use this type of viscosity coefficients.
\begin{definition}
We ask the following properties for the viscosity coefficients $\lambda$ and $\mu$ which are regular:
\begin{itemize}
\item it exists $c>0$, $s_{0}>0$ such that $\forall s$ such that $0\leq s\leq s_{0}$,
$\mu(s)>c$.
\item it exists $c_{1}>0$, $m\in\mathbb{Z}$ such that $\forall s\geq s_{0}$,
$\mu(s)\leq c_{1}s^{m}$.
\item it exists $c^{'}>0$, $s_{0}>0$ such that $\forall s$ such that $0\leq s\leq s^{'}_{0}$,
$\lambda(s)>c$.
\item it exists $c_{2}>0$, $m^{'}\in\mathbb{Z}$ such that $\forall s\geq s^{'}_{0}$,
$\lambda(s)\leq c_{2}s^{m^{'}}$.
\end{itemize}
\label{3defviscosity}
\end{definition}
These hypothesis allow us in fact to control the gradient of the velocity $\n u$ without difficulties. Indeed we have
then:
$$\int_{0}^{+\infty}\int_{\R^{N}}\mu(\rho)|Du|^{2}(t,x)dxdt\geq c\int_{0}^{+\infty}\int_{\R^{N}}|Du|^{2}(t,x)dxdt,$$
and we use the fact that:
$$\int_{0}^{+\infty}\int_{\R^{N}}|Du|^{2}(t,x)dxdt=\frac{1}{2}\big(\int_{0}^{+\infty}\int_{\R^{N}}|\n u|^{2}dxdt
+\int_{0}^{+\infty}\int_{\R^{N}}|{\rm div}u|^{2}dxdt\big),$$
to conclude.
In the following theorem we are interested now by getting a gain of
derivative on the density $\rho$. This will enable us to treat in
distribution sense the quadratic term $\n A(\rho)\otimes\n A(\rho)$.
\begin{theorem}
Let $N=2$ and $(\rho,u)$ be a smooth approximate solution of the system
$(\ref{3systeme})$ with $\kappa(\rho)=\kappa\rho^{\alpha}$ with $\alpha\in\R$ and $\alpha\ne-2$.
We assume that if $\alpha>-2$ then $\frac{1}{\rho}\in
L^{\infty}((0,T)\times\R^{N})$ else $\alpha<-2$ then $\rho\in L^{\infty}((0,T)\times\R^{N})$.\\ Then there exists a constant
$\eta>0$ depending only on the constant coming from the Sobolev
embedding such that if:
$$\|\n\rho_{0}\|_{L^{2}(\R^{2})}+\|\sqrt{\rho_{0}}|u_{0}|\|_{L^{2}(\R^{2})}+\|j_{\gamma}(\rho_{0})\|_{L^{1}}\leq\eta$$
then we get for all $\va\in C^{\infty}_{0}(\R^{N})$:
$$\|\va B(\rho) \|_{L^{2}_{T}(H^{1+\frac{s}{2}})}\leq M\;\;\;\;\;\mbox{with}\;\;0\leq
s<2,$$ where $M$ depends only on the initial conditions data, on
$T$, on $\va$, on $s$ and on $\|\frac{1}{\rho}\|_{L^{\infty}}$ or $\|\rho\|_{L^{\infty}}$.
\label{3T1}
\end{theorem}
\begin{remarka}
When we speak of smooth approximate solution $(\rho,u)$ of the system $(\ref{3systeme})$, in fact we can take
the solution of \cite{3H}. Indeed in this article, we show the existence of global strong solution for the nonisothermal system with capillarity coefficient variable in function of the density. The initial data are in critical Besov spaces for the scaling of the equation.\\
In the sequel when I assume ''smooth approximate solution $(\rho,u)$`` I refer to \cite{3H}. Indeed this is compatible with the fact we will need every time of hypothesis of smallness on the initial data.
\end{remarka}
\begin{remarka}
In fact for the case $\gamma>-2$, we don't need of the hypothesis on the viscosity coefficients concerning
their behavior around the vacuum as we suppose that $\frac{1}{\rho}\in L^{\infty}$. So we can work with degenerate viscosity coefficients.
\end{remarka}
\begin{remarka}
We could remark that in the specific case $\kappa(\rho)=\kappa\rho^{-2}$,  we get a gain of derivative on the density and in particular on $\log\rho$ without condition on the vacuum or on the $L^{\infty}$ norm of the density. So the case $\alpha=-2$ appears as critical.
 In reality the fact that $A(\rho)=B(\rho)$ when $\alpha=-2$ plays a crucial role in the proof. We prove this result in the next theorem.
\end{remarka}
\begin{remarka}
In fact instead of supposing that $\frac{1}{\rho}\in L^{\infty}$ or $\rho\in L^{\infty}$ in theorem \ref{3T1}, we have just to
assume that $\frac{A^{'}(\rho)}{B^{'}(\rho)}\n A(\rho)\in L^{\infty}(L^{2})$.
This would imply that $\rho$ is in fact a weight of Muckenhoupt.
\end{remarka}
\begin{remarka}
In the sequel the notation of space follows those by Runst, Sickel
in \cite{3RS}.
\end{remarka}
{\bf Proof of Theorem \ref{3T1}:}
Our goal is to get a gain of derivative on the density by using
energy inequalities and by taking advantage of the term of
capillarity. We need to localize the argument to control the low
frequencies. Let $\va\in C^{\infty}_{0}(\R^{N})$, we have then by multiplying the momentum equation and applying the operator ${\rm div}$ ( where we use the classical summation index):
\begin{equation}
\begin{aligned}
&\p_{t}{\rm div}(\va\rho u)+\p_{i,j}(\va\rho u_{i}u_{j})-\p_{i,j}(2\va\mu(\rho)Du_{i,j})-\D(\va\lambda(\rho){\rm div}u)+\D(\va P(\rho))\\
&\hspace{1,5cm}=\D\,(A^{1}_{\alpha}\D(\va B(\rho))-A^{2}_{\alpha}\va|\n A(\rho)|^{2})-B_{\alpha}\p^{2}_{i,j}(\va\p_{i}A(\rho)\p_{j}A(\rho))+R_{\va}
\end{aligned}
\label{moment1}
\end{equation}
with:
$$
\begin{aligned}
&R_{\va}=\frac{\p}{\p t}(\rho u\cdot\n\va)+(\p_{i,j}\va)\rho
u_{i}u_{j}+2\p_{i}\va\,\p_{j}(\rho u_{i}u_{j})-(\p_{i,j}\va)\mu(\rho)
Du_{i,j}-A^{1}_{\alpha}\D(B(\rho)\D\va)\\
&\hspace{0,5cm}-2\p_{i}\va\,\p_{j}(\mu(\rho) Du_{i,j})
-\lambda(\rho)
\D\va{\rm div}u-2\n\va\cdot\n(\lambda(\rho){\rm div}u)+\D\va
a\rho^{\gamma}+2a\n\va\cdot\n(\rho^{\gamma})\\
&\hspace{3,1cm}-\D\va(A^{1}_{\alpha}\D B(\rho)-A^{2}_{\alpha}|\n A(\rho)|^{2})
-\n\va\cdot\n(A^{1}_{\alpha}\D B(\rho)-A^{2}_{\alpha}|\n A(\rho)|^{2})\\
&\hspace{2cm}+B_{\alpha}(\p^{2}_{i,j}\va)\p_{i}A(\rho)\p_{j}A(\rho)
+2\kappa\p_{i}\va\,\p_{j}(\p_{i}A(\rho)\p_{j}A(\rho))
-2A^{1}_{\alpha}\D(\n\va\cdot\n B(\rho)).
\end{aligned}
$$
We can apply to the momentum equation the operator
$\Lambda(\D)^{-2}$ in order to make appear a term
in $\Lambda B(\rho)$ coming from the capillarity. Then we obtain:
\begin{equation}
\begin{aligned}
&A^{1}_{\alpha}\Lambda(\va B(\rho))+A^{2}_{\alpha}\Lambda^{-1}(\va|\n A(\rho)|^{2})+
B_{\alpha}\Lambda^{-1}R_{i}R_{j}(\va\p_{i}A(\rho)\p_{j}A(\rho))\\
&=-\Lambda^{-3}\frac{\p}{\p t}{\rm div}(\va\rho\,u)+\Lambda^{-1}R_{i}R_{j}(\va\rho\,u_{i}u_{j})-\Lambda^{-1}(\va\lambda(\rho){\rm
div}u)
+\Lambda^{-1}(\va P(\rho))\\
&\hspace{6cm}-\Lambda^{-1}R_{i}R_{j}(2\mu(\rho)Du_{i,j})+\Lambda^{-1}(\D)^{-1}R_{\va},
\end{aligned}
\label{3C1}
\end{equation}
where $R_{i}$ denotes the classical Riesz operator. We multiply now
the previous equality by $\La^{1+s}(\va B(\rho))$  and we integrate
on space and in time:
\begin{equation}
\begin{aligned}
&A^{1}_{\alpha}\int^{T}_{0}\g|\Lambda^{1+\frac{s}{2}}(\va B(\rho))|^{2}dxdt+A^{2}_{\alpha}\int^{T}_{0}
\g(\va|\n A(\rho)|^{2})\La^{s}(\va B(\rho))dxdt
\\
&+B_{\alpha}\int^{T}_{0}\g\sum_{i,j}R_{i}R_{j}(\va\p_{i}A(\rho)\p_{j}A(\rho))\La^{s}(\va B(\rho))dxdt=\\
&\g\Lambda^{-3}{\rm
div}(\va\rho\,u)\Lambda^{1+s}(\va B(\rho))(T)dx-\g\Lambda^{-3}{\rm
div}(\va\rho_{0}\,u_{0})\Lambda^{1+s}(\va B(\rho_{0}))dx\\
&\hspace{1,3cm}-\int^{T}_{0}\g\big(\Lambda^{-3}{\rm
div}(\va\rho\,u)\La^{1+s}\frac{\p}{\p t}(\va B(\rho))-\va\lambda(\rho){\rm
div}u\,\La^{s}(\va B(\rho))\big)dxdt\\
&\hspace{0,7cm}+\int^{T}_{0}\g\big( R_{i,j}(2\va\mu(\rho)Du_{i,j})\,\La^{s}(\va B(\rho))-R_{i}R_{j}(\va\rho\,u_{i}u_{j})\La^{s}(\va B(\rho))\big)dxdt\\
&\hspace{1,8cm}+\int^{T}_{0}\g \va P(\rho)\La^{s}(\va B(\rho))dxdt
+\int^{T}_{0}\g(\D)^{-1}R_{\va}\La^{s}(\va B(\rho))dxdt.
\label{33M}
\end{aligned}
\end{equation}
Now we want to control the term
$\int^{T}_{0}\g|\Lambda^{1+\frac{s}{2}}(\va B(\rho))|^{2}$. Before
coming in the heart of the proof we want to rewrite the inequality
(\ref{33M}) in particular the term:
$$\int^{T}_{0}\g\Lambda^{-3}{\rm
div}(\va\rho\,u)\La^{1+s}\frac{\p}{\p t}(\va B(\rho)).$$ In this
goal we recall the renormalized equation for $\va B(\rho)$:
\begin{equation}
\frac{\p}{\p t}(\va B(\rho))+{\rm
div}(\va B(\rho)u)=-\va B(\rho){\rm div}u+r_{\va}, \label{3C2}
\end{equation}
with $r_{\va}=-\n\va\cdot B(\rho)u$.
So by using the renormalized equation (\ref{3C2}) we have:
\begin{equation}
\begin{aligned}
&\int^{T}_{0}\g\Lambda^{-3}{\rm div}(\va\rho\,u)\La^{1+s}\frac
{\p}{\p t}(\va B(\rho))dxdt=\\
&-\int^{T}_{0}\g\Lambda^{-2}{\rm
div}(\va\rho\,u)\Lambda^{s}(\va B(\rho){\rm div}u)-\Lambda^{-2}{\rm
div}(\va\rho\,u)\La^{s}{\rm div}(\va B(\rho)u)dxdt\\
&\hspace{7,8cm}+
\int^{T}_{0}\g\Lambda^{-3}{\rm div}(\va\rho\,u)r_{\va}dxdt.
\end{aligned}
\label{3C3}
\end{equation}
By combining (\ref{33M}) and (\ref{3C3}) we get:
$$
\begin{aligned}
&A^{1}_{\alpha}\int^{T}_{0}\g|\Lambda^{1+\frac{s}{2}}(\va B(\rho))|^{2}dxdt+A^{2}_{\alpha}\int^{T}_{0}
\g(\va|\n A(\rho)|^{2})\La^{s}(\va B(\rho))dxdt\hspace{5cm}
\\
&+B_{\alpha}\int^{T}_{0}\g\sum_{i,j}R_{i}R_{j}(\va\p_{i}A(\rho)\p_{j}A(\rho))\La^{s}(\va B(\rho))dxdt=
\end{aligned}
$$
\begin{equation}
\begin{aligned}
&\g\Lambda^{-3}{\rm
div}(\va\rho\,u)\Lambda^{1+s}(\va B(\rho))(T)dx-\g\Lambda^{-3}{\rm
div}(\va\rho_{0}\,u_{0})\Lambda^{1+s}(\va B(\rho_{0}))dx&\\
&-\int^{T}_{0}\g\big(\Lambda^{-2}{\rm
div}(\va\rho\,u)\Lambda^{s}(\va B(\rho){\rm div}u)+\Lambda^{-2}{\rm div}(\va\rho\,u) \La^{s}{\rm
div}(\va B(\rho)u)\big)dxdt&\\
&\hspace{0,5cm}+\int^{T}_{0}\g\big( R_{i,j}(2\va\mu(\rho)Du_{i,j})\,\La^{s}(\va B(\rho))-R_{i}R_{j}(\va\rho\,u_{i}u_{j})\La^{s}(\va B(\rho))\big)dxdt&\\
&\hspace{1,2cm}+\int^{T}_{0}\g\va\lambda(\rho){\rm
div}u\,\La^{s}(\va B(\rho))dxdt+\int^{T}_{0}\g
\va P(\rho)\La^{s}(\va B(\rho))dxdt&\\
&\hspace{0,5cm}+\int^{T}_{0}\g(\D)^{-1}R_{\va}\La^{1+\frac{s}{2}}(\va B(\rho))dxdt
+\int^{T}_{0}\g\La^{-3}{\rm div}(\va\rho u)\La^{1+s}r_{\va}dxdt.&
\label{3a3}
\end{aligned}
\end{equation}
In order to control
$\int^{T}_{0}\g|\Lambda^{1+\frac{s}{2}}(\va B(\rho))|^{2}$, it
suffices to bound all the other terms of (\ref{3a3}). It will allow us to get
a control on $\La^{1+\frac{s}{2}}(\va B(\rho))$ and so a gain of $\frac{s}{2}$ derivative on the gradient of density $\n A(\rho)$.\\
We start with the most complicated term which requires a control
of $\frac{1}{\rho}$ in $L^{\infty}$ if $\alpha>-2$ and a control of $\rho\in L^{\infty}$ if $\alpha<-2$.
In the sequel, we will treat only the case $\alpha>-2$, the proof of the other case follows exactly the same lines.
So the fact that $\frac{1}{\rho}$ is in $L^{\infty}$ will imply in particular that $\n u\in L^{2}(L^{2})$ and $u\in L^{\infty}(L^{2})$, this fact will be permanently use in the sequel.
\subsubsection*{1)\;$\int_{0}^{T}\g(\va|\n A(\rho)|^{2})\La^{s}(\va B(\rho))$:}
By induction we have $\n(\va B(\rho))\in
L^{2}_{T}(\dot{H}^{\frac{s}{2}})$
%
%
and by using Sobolev embedding we get $\n(\va B(\rho))\in L^{2}(L^{p})$
with
$\frac{1}{p}=\frac{1}{2}-\frac{s}{4}$ (we remark that the case $s=2$ is critical for Sobolev embedding).
Now we have:
$$\va\n A(\rho)=\frac{(\frac{\alpha}{2}+1)\n(\va B(\rho))}{(2+\alpha)A(\rho)}-
\frac{(\frac{\alpha}{2}+1)A(\rho)}{(2+\alpha)}\n\va\;\;\;\;\mbox{if}\;\;\alpha\ne-2.$$
and
we recall that by hypothesis $\frac{1}{\rho}\in L^{\infty}$, so we
have $\va\n A(\rho)\in L^{2}(L^{p})$ because $\frac{(\frac{\alpha}{2}+1)B(\rho)}{(2+\alpha)A(\rho)}\n\va\in
L^{\infty}(L^{r})$
for all $1\leq r\leq+\infty$ as $\n A(\rho)\in L^{\infty}(L^{2})$.\\
We now consider $\La^{s}(\va B(\rho))$. We have by induction
$\La^{s}(\va B(\rho))\in L^{2}(\dot{H}^{1-\frac{s}{2}})$ and
$\La^{s}(\va B(\rho))$ is in $ L^{2}(L^{2})$ because $\va B(\rho)\in
L^{2}(L^{2})$ (here the fact to localize by $\va$ is crucial) which enables us to control the low frequencies of
$\La^{s}(\va B(\rho))$. We have then $\La^{s}(\va B(\rho))\in
L^{2}(H^{1-\frac{s}{2}})\h L^{2}(L^{q})$ with $\frac{1}{q}=\frac{s}{4}$.
Finally by H\"older inequality we get
$\va|\n A(\rho)|^{2}\La^{s}(\va B(\rho))\in L^{1}_{T}(L^{1}(\R^{N}))$
because
$\frac{1}{2}+\frac{1}{p}+\frac{1}{q}=\frac{1}{2}+\frac{1}{2}-\frac{s}{4}+\frac{s}{4}=1$
and we get more precisely:
\begin{equation}
\begin{aligned}
&\int_{0}^{T}\g(\va|\n A(\rho)|^{2})\La^{s}(\va B(\rho))dxdt\lesssim\|\n A(\rho)\|_{L^{\infty}_{T}(L^{2})}
\|\La^{s}(\va B(\rho))\|_{L^{2}_{T}(L^{q})}\\
&\hspace{9,5cm}\times\|\va\n A(\rho)\||_{L^{2}_{T}(L^{p})},\\
&\hspace{0cm}\lesssim\|\frac{1}{\rho}\|_{L^{\infty}_{T}(L^{\infty})}
\|\Lambda^{1+\frac{s}{2}}(\va B(\rho))\|_{L^{2}_{T}(L^{2})}\|\n A(\rho)\|_{L^{\infty}_{T}(L^{2})}
(1+\|\Lambda^{1+\frac{s}{2}}(\va B(\rho))\|_{L^{2}_{T}(L^{2})}).
\end{aligned}
\end{equation}
We proceed similarly for the term:
$$\int_{0}^{T}\g\sum_{i,j}R_{i}R_{j}(\va\p_{i}A(\rho)\p_{j}A(\rho))\La^{s}(\va B(\rho))dxdt,$$
indeed we have in following the same lines
$\va\p_{i}A(\rho)\p_{j}A(\rho)\in L^{2}(L^{q})$ with
$\frac{1}{q}=1-\frac{s}{4}$ and we use the fact that the Riesz operator is
continuous from $L^{p}$ in $L^{p}$ for $1<p<+\infty$.
\\
We next study the term $\g\La^{-3}{\rm
div}(\va\rho\,u)\La^{1+s}(\va B(\rho))(t)dxdt$.
\subsubsection*{2) $\g\La^{-3}{\rm div}(\va\rho\,u)\La^{1+s}(\va B(\rho))dx$:}
We rewrite the term $\g\La^{-3}{\rm
div}(\va\rho\,u)\La^{1+s}(\va B(\rho))$ on the following form:
$$
\begin{aligned}
\g\La^{-3}{\rm
div}(\va\rho\,u)\La^{1+s}(\va B(\rho))dx=\g&\La^{-1}{\rm
div}(\va\rho\,u)\La^{-1+s}(\va B(\rho))dx\\
&=\sum_{1\leq i\leq N}\g R_{i}(\va\rho
u_{i})\La^{-1+s}(\va B(\rho))dx.
\end{aligned}
$$
As $\frac{1}{\rho}\in L_{T}^{\infty}(L^{\infty})$ then we have $u\in
L^{\infty}_{T}(L^{2})$. We recall that $\n A(\rho)\in
L^{\infty}(H^{1})$ then $\va\rho\in L^{\infty}(L^{p})$ for all
$1\leq p<+\infty$. We deduce that $\va\rho u$ belongs to
$L^{\infty}(L^{2-\beta}\cap L^{1})$ for $\beta>0$.
So we have $R_{i}(\va\rho u_{i})\in L^{\infty}_{T}(L^{r})$ for all $1<r<2$ by continuity of the operator $R_{i}$ from $L^{p}$ to
$L^{p}$ when $1<p<+\infty$.
\subsubsection*{Case $1\leq s<2$:}
Next we have:
$$\n(\va B(\rho))=\va B^{'}(\rho)\n\rho+B(\rho)\n\va$$
then we get $\n(\va B(\rho))\in L^{\infty}(L^{2-\beta})$, by using
the fact that $\n A(\rho)\in L^{\infty}(L^{2})$ and Sobolev
embedding with H\"older inequalities.
We have then that $\va B(\rho)$ belongs to $L^{\infty}(W^{1}_{2-\beta})$.
So $\La^{s-1}(\va B(\rho))$ belongs to $
L^{\infty}(H^{2-s}_{2-\beta})$. By Sobolev embedding
$\La^{s-1}(\va B(\rho))$ is in  $L^{\infty}(L^{p})$ with
$\frac{1}{p}=\frac{1}{2-\beta}-\frac{2-s}{2}=\frac{1}{2-\beta}+\frac{s}{2}-1$
with $\beta$ small enough to avoid critical
embedding.
Finally we get $R_{i}(\va\rho u)\La^{-1+s}(\va B(\rho))\in
L^{1}_{T}(L^{1}(\R^{N})$. Indeed we have
$\frac{1}{p}+\frac{1}{2-\beta}=\frac{2}{2-\beta}+\frac{s}{2}-1<1$ by
choosing $\beta$ small enough and
$\frac{1}{p}+\frac{1}{1+\beta}=\frac{1}{1+\beta}-1+\frac{2}{2-\beta}+\frac{s}{2}>1$
by choosing $\beta$ small big if necessary, we conclude by
interpolation. Here the fact that $\va$ is in a compact support is crucial.
We have finally:
$$\big|\g\La^{-3}{\rm div}(\va\rho\,u)\La^{1+s}(\va B(\rho))dx\big|\leq
M_{0},$$ with $M_{0}$ depending only on the initial data.
\subsubsection*{Case $0<s<1$:}
In this case we conclude by interpolation with the previous case.
We now want to study the other terms coming from the renormalized
equation (\ref{3C2}).
\subsubsection*{3) $\int_{0}^{T}\g\La^{-2}{\rm div}(\va\rho\,u)\La^{s}(\va B(\rho){\rm div}u)dxdt$,
$\int_{0}^{T}\g\La^{-2}{\rm div}(\va\rho\,u)\La^{s}({\rm
div}(\va B(\rho)\,u))dxdt$:}
We start with:
$$\int_{0}^{T}\g\La^{-1}{\rm div}(\va\rho\,u)\La^{s-1}({\rm div}(\va B(\rho)\,u))=
\int_{0}^{T}\g{\rm div}(\va\rho\,u)\La^{s-2}({\rm
div}(\va B(\rho)\,u)).$$
\subsubsection*{Case $1\leq s<2$:}
We have:
$${\rm div}(\va\rho\,u)=u\cdot\n(\va\rho)+\va\rho{\rm div}u.$$
By H\"older inequalities and Sobolev embedding we get that ${\rm
div}(\va\rho\,u)$ belongs to $L^{2}_{T}(L^{2-\beta})$ for all
$\beta\in]0,1]$.
Next we rewrite ${\rm div}(\va B(\rho)\,u)$ on the form:
$${\rm div}(\va B(\rho)\,u)=u\cdot\n(\va B(\rho))+\va B(\rho){\rm
div}u.$$
As previously ${\rm div}(\va B(\rho)\,u)$ is in
$L^{2}_{T}(L^{2-\beta})$ for all $\beta\in]0,1]$. Now by Sobolev
embedding we have $\La^{s-2}{\rm div}(\va B(\rho)\,u)\in
L^{2}_{T}(L^{p})$ with $\frac{1}{p}=\frac{1}{2-\beta}-\frac{2-s}{2}$
with $\beta$ small enough to avoid critical Sobolev embedding.
We conclude that ${\rm div}(\va\rho\,u)\La^{s-2}({\rm
div}(\va B(\rho)\,u))$ is in $L^{1}_{T}(L^{1})$ because
$\frac{1}{2-\beta}+\frac{1}{p}=\frac{2}{2-\beta}-\frac{2-s}{2}=\frac{2}{2-\beta}-1+\frac{s}{2}<1$
with $\beta$ small enough if necessary and $1+\frac{1}{p}>1$, so we
obtain the result by interpolation. Finally we have:
$$\big|\int_{0}^{T}\g\La^{-2}{\rm div}(\va\rho\,u)\La^{s}({\rm
div}(\va B(\rho)\,u))dxdt\big|\leq M_{0},$$ with $M_{0}$ depending
only on the initial data.
\subsubsection*{Case $0<s<1$:}
We have the result by interpolation with the previous case.
Next we proceed similarly for:
$$\int_{0}^{T}\g\La^{-1}{\rm div}(\va\rho\,u)\La^{s-1}(\va B(\rho){\rm div}u)dxdt.$$
\subsubsection*{4) Last terms}
We now want to concentrate us on the following term:
$$\int^{T}_{0}\g\sum_{i,j}R_{i}R_{j}(\va\rho\,u_{i}u_{j})\La^{s}(\va B(\rho))dxdt.$$
We know that $u\in L^{\infty}(L^{2})$ (as $\frac{1}{\rho}\in L^{\infty}$) and $D u\in L^{2}(L^{2})$ then
$u\in L^{2}_{T}(H^{1})$ and by H\"older inequalities and Sobolev
embedding we can show that $\va\rho\,u_{i}u_{j}\in L^{2}_{T}(L^{2-\beta})$ with $\beta>0$
and so $R_{i}R_{j}(\va\rho\,u_{i}u_{j})\in L^{2}_{T}(L^{2-\beta})$.\\
We have seen that $\La^{s}(\va B(\rho))\in L^{2}_{T}(H^{1-\frac{s}{2}})$
then we have as $1-\frac{s}{2}>0$ for $\beta$ small enough:
$$\|\La^{s}(\va B(\rho))\|_{L^{2}_{T}(L^{\frac{2-\beta}{1-\beta}})}\leq M_{0}+\|\va B(\rho)\|_{L^{2}_{T}(\dot{H}^{1-\frac{s}{2}})}^{\alpha}$$
with $0<\alpha<1$.
We have then:
$$|\int^{T}_{0}\g\sum_{i,j}R_{i}R_{j}(\va\rho\,u_{i}u_{j})\La^{s}(\va B(\rho))dxdt|\leq M_{0}
+\|\La^{1+\frac{s}{2}}(\va
B(\rho))\|^{\beta} _{L^{2}(L^{2})}$$
with $0<\beta<1$ and $M_{0}$ depending only on the initial data.\\
We are interested in the term:
$\int_{0}^{t}\g\va\lambda(\rho){\rm div}u\,\La^{s}(\va B(\rho))dxdt$
We have then ${\rm div}u\in L^{2}(L^{2})$ and we have shown that
$\La^{s}(\va B(\rho))\in L^{2}_{T}(H^{1-\frac{s}{2}})$ so we conclude in the
same way than the previous term.\\ We proceed similarly for the term:
$\int^{T}_{0}\g R_{i,j}(2\va\mu(\rho)Du_{i,j})\,\La^{s}(\va B(\rho))dxdt$.
We finally conclude with the term:
$$\int_{0}^{T}\g \va P(\rho)\La^{s}(\va B(\rho))dxdt.$$
Similarly we have $\La^{s}(\va B(\rho))\in L^{2}_{T}(L^{2})$ and
$\va P(\rho)\in L^{2}_{T}(L^{2})$ because $\va A(\rho)$ is in
$L^{\infty}(H^{1})$, and we conclude by Sobolev embedding.
To finish we have to control the term depending on $R_{\va}$ and $r_{\va}$
that we leave to the reader. Indeed these terms are easy to treat because they are
more regular than the previous terms.\\
We finally get by using all the previous inequalities:
$$
\begin{aligned}
&\|\La^{1+\frac{s}{2}}(\va B(\rho))\|^{2}_{L^{2}(L^{2})}\leq C_{0}\|\n A(\rho)\|_{L^{\infty}(L^{2})}\|\La^{1+
\frac{s}{2}}(\va B(\rho))\|^{2}_{L^{2}(L^{2})}\\
&\hspace{8cm}+C_{1}\|\La^{1+\frac{s}{2}}(\va B(\rho))\|^{2\beta}_{L^{2}(L^{2})}+M_{0}
\end{aligned}
$$
with $0<\beta<1$ and $ C_{0}$, $C_{1}$, $M_{0}$ depends only of the norm of initial data.
By energy inequalities we have
$\|\n A(\rho)\|_{L^{\infty}(L^{2})}\leq\e<1$, we can then conclude that:
$$\|\va B(\rho)\|_{L^{2}(\dot{H}^{1+\frac{s}{2}})}\leq M_{0}$$
with $M_{0}$ depending only on the initial data and $\va$. \hfill {$\Box$}
\\
We want now improve this result by extracting a specific structure of the capillarity tensor.
Indeed by choosing $\kappa(\rho)=\kappa\rho^{-2}$ with $\kappa>0$, we are going show that we get the same estimate but without any conditions on the vacuum. In fact, the power $-2$ seems to be critical, and it explains very easily in our proof by the simple fact that $A(\rho)=B(\rho)$. So we can apply a boobstrap argument without asking a control of $\rho$ or $\frac{1}{\rho}$ in $L^{\infty}$.\\
However a new difficulty appears on the control of the pressure term as in our case we have only a control
on $\n\ln\rho\in L^{\infty}(L^{2})$. So we can hope directly a control of $\va\rho\in L^{\infty}(L^{p})$ for all $1\leq p<+\infty$.
We will see in fact that we can in the same time as a gain of derivative on the density obtaining a gain of integrability on the density.
\begin{theorem}
Let $N=2$, $\kappa>0$ and $(\rho,u)$ be a smooth approximate solution of the system
$(\ref{3systeme})$ with $\kappa(\rho)=\kappa\rho^{-2}$.Then there exists a constant
$\eta>0$ depending only on the constant coming from the Sobolev
embedding such that if:
$$\|\n\rho_{0}\|_{L^{2}(\R^{2})}+\|\sqrt{\rho_{0}}|u_{0}|\|_{L^{2}(\R^{2})}+\|j_{\gamma}(\rho_{0})\|_{L^{1}}\leq\eta$$
then it exists $\alpha>0$ such that for all $\va\in C^{\infty}_{0}(\R^{N})$:
$$\|\va B(\rho) \|^{2}_{L^{2}_{T}(H^{1+\frac{s}{2}})}+\|\va\rho^{\alpha-2}\n\rho\|^{2}_{L^{2}(L^{2})}+\|\va\rho\|_{L^{\gamma+\alpha}
((0,T)\times\R^{N})}^{\gamma+\alpha}\leq M\;\;\;\;\;\mbox{with}\;\;0\leq
s\leq \e,$$ where $M$ depends only on the initial conditions data, on
$T$, on $\va$ and on $\e$. $\e$ depends only of $\gamma$ the coefficient of the pressure and is small.
\label{a3T1}
\end{theorem}
{\bf Proof:} The proof follow the same line as the proof of theorem \ref{3T1} except concerning the bounds of
estimate coming from the capillarity term and the fact that we lost the control of $\rho$ in $L^{\infty}(L^{p}_{loc})$ for all $1\leq p<+\infty$.
We need then to get a gain of integrability on the density to treat the term coming of the pressure.
We apply to equation (\ref{moment1}) the operator $\Lambda$ and we multiply par $\rho^{\alpha}$ with $\alpha>0$:
\begin{equation}
\begin{aligned}
&\big[A^{1}_{\alpha}\D(\va B(\rho))+A^{2}_{\alpha}(\va|\n A(\rho)|^{2})+
B_{\alpha}R_{i}R_{j}(\va\p_{i}A(\rho)\p_{j}A(\rho))\big]\rho^{\alpha}\\
&=\frac{\p}{\p t}(\D)^{-1}{\rm div}(\va\rho\,u)\rho^{\alpha}+R_{i}R_{j}(\va\rho\,u_{i}u_{j})\rho^{\alpha}-\va\lambda(\rho)({\rm
div}u)\rho^{\alpha}
+\va\rho^{\gamma+\alpha}\\
&\hspace{6cm}-R_{i}R_{j}(2\mu(\rho)Du_{i,j})\rho^{\alpha}+(\D)^{-1}R_{\va}\rho^{\alpha},
\end{aligned}
\label{3C11}
\end{equation}
Next we integrate on $(0,T)\times\R^{N}$ the equation \ref{3C11} we get then:
\begin{equation}
\begin{aligned}
&\int_{(0,T)\times\R^{N}}\va\rho^{\gamma+\alpha}(x,t)dxdt+A_{\alpha}\int_{(0,T)\times\R^{N}}\rho^{\alpha-2}
|\n\rho|^{2}dxdt=\\
&
B_{\alpha}\int_{(0,T)\times\R^{N}}\va\p_{i}\ln(\rho)\p_{j}\ln(\rho)R_{i,j}\rho^{\alpha}dxdt
-\int_{\R^{N}}(\D)^{-1}{\rm div}(\va\rho\,u)\rho^{\alpha}(T,x)dx\\
&+\int_{\R^{N}}(\D)^{-1}{\rm div}(\va\rho_{0}\,u_{0})\rho_{0}^{\alpha}dx-\int_{(0,T)\times\R^{N}}(\D)^{-1}{\rm div}(\va\rho\,u)\p_{t}\rho^{\alpha}dxdt\\
&-\int_{(0,T)\times\R^{N}}R_{i,j}(\va\rho\,u_{i}u_{j})\rho^{\alpha}dxdt-\int_{(0,T)\times\R^{N}}\va\lambda(\rho)({\rm
div}u)\rho^{\alpha}dxdt
\\
&-\int_{(0,T)\times\R^{N}}R_{i}R_{j}(2\mu(\rho)Du_{i,j})\rho^{\alpha}dxdt+\int_{(0,T)\times\R^{N}}(\D)^{-1}R_{\va}\rho^{\alpha}dxdt.
\end{aligned}
\label{C4}
\end{equation}
In the sequel we will note:
$$
\begin{aligned}
&F=-\int_{\R^{N}}(\D)^{-1}{\rm div}(\va\rho\,u)\rho^{\alpha}(T,x)dx
+\int_{\R^{N}}(\D)^{-1}{\rm div}(\va\rho_{0}\,u_{0})\rho_{0}^{\alpha}dx\\
&-\int_{(0,T)\times\R^{N}}(\D)^{-1}{\rm div}(\va\rho\,u)\p_{t}\rho^{\alpha}dxdt-\int_{(0,T)\times\R^{N}}R_{i,j}(\va\rho\,u_{i}u_{j})\rho^{\alpha}dxdt\\
&-\int_{(0,T)\times\R^{N}}\va\lambda(\rho)({\rm
div}u)\rho^{\alpha}dxdt
-\int_{(0,T)\times\R^{N}}R_{i}R_{j}(2\mu(\rho)Du_{i,j})\rho^{\alpha}dxdt\\
&+\int_{(0,T)\times\R^{N}}(\D)^{-1}R_{\va}\rho^{\alpha}dxdt.
\end{aligned}
$$
Our goal is now to have $\rho\in L^{\gamma+\alpha}_{loc}((0,T)\times\R^{N})$, for this we have to control all the terms on the right handside. This processus follow the same lines that this one used in the case of Navier-Stokes compressible to get a gain of integrability on the pressure. The new difficulty consists only in the following term:
$$\int_{(0,T)\times\R^{N}}\va\p_{i}\ln(\rho)\p_{j}\ln(\rho)R_{i,j}\rho^{\alpha}dxdt,$$
To control this term we need of a gain of derivative on $\n\ln(\rho)$ that is why we have then in suming equation (\ref{C4})
and (\ref{33M}):
$$
\begin{aligned}
&A^{1}_{\alpha}\int^{T}_{0}\g|\Lambda^{1+\frac{s}{2}}(\va B(\rho))|^{2}dxdt+\int_{(0,T)\times\R^{N}}\big(\va\rho^{\gamma+\alpha}(t,x)+A_{\alpha}\rho^{\alpha-2}|\n\rho|^{2}\big)dxdt=\\
&-A^{2}_{\alpha}\int^{T}_{0}
\g(\va|\n A(\rho)|^{2})\La^{s}(\va B(\rho))dxdt-B_{\alpha}\int_{(0,T)\times\R^{N}}\va\p_{i}\ln(\rho)\p_{j}\ln(\rho)R_{i,j}\rho^{\alpha}dxdt\\
&-B_{\alpha}\int^{T}_{0}\g\sum_{i,j}R_{i}R_{j}(\va\p_{i}A(\rho)\p_{j}A(\rho))\La^{s}(\va B(\rho))dxdt+F\\
&-\g\Lambda^{-3}{\rm
div}(\va\rho_{0}\,u_{0})\Lambda^{1+s}(\va B(\rho_{0}))dx+\g\Lambda^{-3}{\rm
div}(\va\rho\,u)\Lambda^{1+s}(\va B(\rho))(T)dx\\
&-\int^{T}_{0}\g\big(\Lambda^{-3}{\rm
div}(\va\rho\,u)\La^{1+s}\frac{\p}{\p t}(\va B(\rho))-
\va\lambda(\rho){\rm
div}u\,\La^{s}(\va B(\rho))\big)dxdt\\
&+\int^{T}_{0}\g\big( R_{i,j}(2\va\mu(\rho)Du_{i,j})\,\La^{s}(\va B(\rho))-R_{i}R_{j}(\va\rho\,u_{i}u_{j})\La^{s}(\va B(\rho))\big)dxdt\\
&+\int^{T}_{0}\g \va P(\rho)\La^{s}(\va B(\rho))dxdt+\int^{T}_{0}\g(\D)^{-1}R_{\va}\La^{s}(\va B(\rho))dxdt.
\label{33}
\end{aligned}
$$
We have now to control all the term on the right handside.
We begin with the capillarity term:
$$\int^{T}_{0}\va|\n(\ln\rho)|^{2}\Lambda^{s}(\va\ln(\rho))dxdt,$$
We have then by induction and Sobolev embedding $\n(\va\ln\rho)\in L^{2}(L^{p})$ and $\Lambda^{s}(\va\ln(\rho))$ is in $L^{2}(L^{p_{1}})$
where $\frac{1}{p}=\frac{1}{2}-\frac{s}{4}$ and  $\frac{1}{p_{1}}=\frac{s}{4}$
so by H\"older inequalities we have a control of $\va|\n(\ln\rho)|^{2}\Lambda^{s}(\va\ln(\rho))$ in $L^{1}_{T}(L^{1})$ because
$\frac{1}{2}+\frac{s}{4}+\frac{1}{2}-\frac{s}{4}=1$.
Next we treat the following term:
$$\int_{(0,T)\times\R^{N}}\va\p_{i}\ln(\rho)\p_{j}\ln(\rho)R_{i,j}(\rho^{\alpha})dxdt,$$
We have then by induction $R_{i,j}(\va\rho^{\alpha})\in L^{\frac{\gamma+\alpha}{\alpha}}((0,T)\times\R^{N})$ and $\va\p_{i}(\ln\rho)\in L^{2}(L^{p})$
with $p=\frac{4}{s}$, by H\"older inequalities we have $\va\p_{i}\ln(\rho)\p_{j}\ln(\rho)R_{i,j}\rho^{\alpha}\in L^{1}(L^{1})$ because
$\frac{1}{2}+\frac{s}{4}+\frac{\alpha}{\gamma+\alpha}\leq1$ and $\frac{1}{2}+\frac{\alpha}{\gamma+\alpha}\leq1$.\\
Similarly we have $\va P(\rho)\in L^{\frac{\gamma+\alpha}{\gamma}}((0,T)\times\R^{N})$ and $\Lambda^{s}\ln(\va\rho)\in L^{\infty}(L^{\frac{2}{s}})$ by Sobolev embedding, we conclude by H\"older inequalities as $\frac{\gamma}{\alpha+\gamma}+\frac{s}{2}\leq 1$.
To finish we 
study the term asking control $L^{\infty}$ in time as
$\g\Lambda^{-1}{\rm
div}(\va\rho\,u)\Lambda^{-1+s}(\va \ln(\rho))(T)dx$ coming from $F$, we have: $\Lambda^{-1}{\rm
div}(\va\rho\,u)\in L^{\infty}(L^{p})$ and $\Lambda^{1+s}(\va \ln(\rho))\in L^{\infty}(L^{q})$
with $\frac{1}{p}=\frac{1}{2\gamma}+\frac{1}{2}$ and $\frac{1}{q}=\frac{1}{2}-\frac{2-s}{2}=\frac{s}{2}-\frac{1}{2}$, we have then
$\Lambda^{-1}{\rm
div}(\va\rho\,u)\Lambda^{-1+s}(\va \ln(\rho))(T)\in L^{\infty}(L^{1})$ if $\frac{1}{2\gamma}+\frac{s}{2}\leq 1$.
The others terms are left to the reader.
We conclude by writing the final estimate where we set:
$$
\begin{aligned}
&A(T)=\int_{(0,T)\times\R^{N}}\va\rho^{\gamma+\alpha}(x,t)dxdt+A_{\alpha}\int_{(0,T)\times\R^{N}}\rho^{\alpha-2}|\n\rho|^{2}dxdt,\\
&B(T)=\int^{T}_{0}\g|\Lambda^{1+\frac{s}{2}}(\va B(\rho))|^{2}dxdt.
\end{aligned}
$$
We have then the following final estimate:
$$A(T)+B(T)\leq\e B(T)+\e A(T)^{\frac{\alpha}{\alpha+\gamma}}B(T)^{\frac{1}{2}}+C A(T)^{\beta_{1}}+C^{'} B(T)^{\beta_{2}}+M,$$
with $0<\beta_{1},\beta_{2}<1$. By boosbstrap we can conclude.
\hfill {$\Box$}
\begin{remarka}
The major difficulty in the previous proof is to treat the terms coming from the capillarity, it implies to
impose a smallness condition on the initial data.
An other idea would be to use refined Sobolev embedding to avoid smallness condition. So we have:
$$\int_{(0,T)\times\R^{N}}\va|\n\ln(\rho)|^{2}\Lambda^{s}\ln(\rho)dxdt\leq\|\va\n\ln(\rho)\|_{L^{2}(L^{p_{1}})}
\|\Lambda^{s}\ln(\rho)\|_{L^{2}(L^{q_{1}})}\|\n\ln(\rho)\|_{L^{\infty}(L^{2})},$$
with $\frac{1}{p_{1}}=\frac{1}{2}-\frac{s}{4}$ and $\frac{1}{q_{1}}=\frac{s}{4}$. Next by using propositions \ref{3refined}
we get:
$$
\begin{aligned}
&\|\va\n\ln(\rho)\|_{L^{p}}\leq C\|\va\ln\rho\|^{\frac{s}{2}}_{B^{\frac{s}{2}}_{\infty,\infty}}\|\va\ln\rho\|^{\frac{2-s}{2}}_{H^{1+\frac{s}{2}}}\\
&\|\va\Lambda^{s}\ln(\rho)\|_{L^{p}}\leq C\|\va\ln\rho\|^{\frac{2-s}{2}}_{B^{\frac{s}{2}}_{\infty,\infty}}\|\va\ln\rho\|^{\frac{s}{2}}_{H^{1+\frac{s}{2}}}
\end{aligned}
$$
We have then:
$$\|\va\n\ln(\rho)\|_{L^{p_{1}}}\|\va\Lambda^{s}\ln(\rho)\|_{L^{q_{1}}}\leq C^{2}
\|\va\ln\rho\|_{B^{\frac{s}{2}}_{\infty,\infty}}\|\va\ln\rho\|_{H^{1+\frac{s}{2}}}$$
We could conclude our argument if $\va\ln\rho\in L^{2}(B^{\e}_{\infty,\infty})$ with $\e$ arbitraly small. In this case
we would be able to get a gain of derivative on the density without any condition of smallnesse on the initial data.
We can notice that in \cite{3DD}, $\rho$ is in $L^{2}(B^{1}_{\infty,1})$ so it's wide enough.
\end{remarka}
In the following result, we want show that we can obtain similar result when we approximate the capillarity coefficient $\kappa(\rho)$ by a constant $\kappa$.
\begin{corollary1}
\label{corollaire2}
Let $N=2$ and $\alpha,M,\kappa\in\R$. $(\rho,u)$ is a smooth approximate solution of the system
$(\ref{3systeme})$ with the following capillarity coefficient:
$$\kappa(\rho)=\frac{1}{\rho^{2}}1_{\{\rho<\alpha\}}+\theta_{1}(\rho)1_{\{\alpha\leq\rho\leq2
\alpha\}}+\kappa 1_{\{2\alpha<\rho\}}.$$
where $\theta_{1},\;\theta_{2}$ are regular function such that $\kappa$ is a regular function.
Then there exists a constant
$\eta>0$ depending only on the constant coming from the Sobolev
embedding such that if:
$$\|\n\rho_{0}\|_{L^{2}(\R^{2})}+\|\sqrt{\rho_{0}}|u_{0}|\|_{L^{2}(\R^{2})}+\|j_{\gamma}(\rho_{0})\|_{L^{1}}\leq\eta$$
then we get for all $\va\in C^{\infty}_{0}(\R^{N})$:
$$\|\va B(\rho)\|_{L^{2}_{T}(\dot{H}^{1+\frac{s}{2}})}\leq M\;\;\;\;\mbox{with}\;\;0<s<2,$$
where $M$ depends only on the initial conditions data, on
$T$, on $\va$ and on $s$.
\end{corollary1}
{\bf Proof:}
The proof follows the same lines than proof of theorem \ref{3T1}. We have then in this case:
$$
\begin{aligned}
&A(\rho)=\log\rho 1_{\{\rho<\alpha\}}+\beta_{1}(\rho)1_{\{\alpha\leq\rho\leq2
\alpha\}}+\kappa \rho 1_{\{2\alpha<\rho\}},\\
&B(\rho)=\log\rho 1_{\{\rho<\alpha\}}+\beta_{2}(\rho)1_{\{\alpha\leq\rho\leq2
\alpha\}}+\rho^{2}1_{\{2\alpha<\rho\}}.
\end{aligned}
$$
where $\beta_{1}(x)=\int^{x}_{0}\sqrt{\theta_{1}}(y)1_{\{\alpha\leq\rho\leq2
\alpha\}}dy$ and $\beta_{2}(x)=\int^{x}_{0}y\theta_{1}(y)1_{\{\alpha\leq\rho\leq2
\alpha\}}dy$.
 Here to simplify we assume that the viscosity coefficient are constant. At the difference with the previous proof, this case is more simple because we have $\n(\rho 1_{\{2\alpha<\rho\}})\in L^{\infty}(L^{2})$, so we have easily $\rho\in L^{\infty}(L^{p}_{loc})$
for all $1\leq p<+\infty$. It implies that we don't need of condition of smallness on s as in the previous case.\\
We apply now to the equation (\ref{moment1}) the operator
$(\D)^{-1}{\rm div}$, next we multiply by $\D(B(\rho))$
and we integrate on space and in time:
$$
\begin{aligned}
&\int^{T}_{0}\g\big(|\Lambda^{1+\frac{s}{2}}B(\rho)|
^{2}+|\n A(\rho)||\frac{A^{'}(\rho)}{B^{'}(\rho)}(\n(\va B(\rho))-B(\rho)\n\va)|
\Lambda^{s}(\va B(\rho)\big)dxdt
\\
&+\int^{T}_{0}\g\sum_{i,j}R_{i,j}(\kappa(\rho)\va\p_{i}\rho\p_{j}\rho)\La^{s}(\va B(\rho))dxdt=
\g\Lambda^{-3}{\rm
div}(\va\rho\,u)\Lambda^{1+s}(\va B(\rho))(T)dx\\
&\hspace{0,5cm}-\g\Lambda^{-3}{\rm
div}(\va\rho_{0}\,u_{0})\Lambda^{1+s}(\va B(\rho_{0}))dx-\int^{T}_{0}\g\Lambda^{-3}{\rm
div}(\va\rho\,u)\La^{1+s}\frac{\p}{\p t}(\va B(\rho))dxdt\\
&+(2\mu+\lambda)\int^{T}_{0}\g\va{\rm
div}u\,\La^{s}(\va B(\rho))dxdt-\int^{T}_{0}\g\sum_{i,j}R_{i}R_{j}(\va\rho\,u_{i}u_{j})\La^{s}(\va B(\rho))dxdt\\
&\hspace{2,9cm}+\int^{T}_{0}\g \va P(\rho)\La^{s}(\va B(\rho))dxdt
+\int^{T}_{0}\g(\D)^{-1}R_{\va}\La^{s}(\va B(\rho))dxdt.
\label{33}
\end{aligned}
$$
Now we proceed as in the proof of theorem \ref{3T1} and we have to control the terms:
$$\int^{T}_{0}
\g(|\n A(\rho)|(\frac{A^{'}(\rho)}{B^{'}(\rho)}(\n(\va B(\rho))-B(\rho)\n\va)
\Lambda^{s}(\va B(\rho))dxdt,$$
Tre crucial point is that $\frac{A^{'}(\rho)}{B^{'}(\rho)}$ is in $L^{\infty}(L^{\infty})$. So as previously we have:
$$
\begin{aligned}
&\int^{T}_{0}
\g(|\n A(\rho)|\frac{A^{'}(\rho)}{B^{'}(\rho)}(\n(\va B(\rho))
\Lambda^{s}(\va B(\rho))dxdt\\
&\hspace{0,7cm}\leq C \|\frac{A^{'}(\rho)}{B^{'}(\rho)}\|_{L^{\infty}(L^{\infty})}\|\n A(\rho)\|_{L^{\infty}(L^{2})}\|\n(\va B(\rho)\|_{L^{2}(L^{p})}\|\Lambda^{s}(\va B(\rho)\|_{L^{2}(L^{q})}, \\
&\hspace{0,7cm}\leq C \|\frac{A^{'}(\rho)}{B^{'}(\rho)}\|_{L^{\infty}(L^{\infty})}\|\n A(\rho)\|_{L^{\infty}(L^{2})}\|\Lambda^{1+\frac{s}{2}}(\va B(\rho)\|^{2}_{L^{2}(L^{2})},
\end{aligned}
$$
with $\frac{1}{p}=1-\frac{s}{4}$ and $\frac{1}{q}=\frac{s}{4}$.
We get finally:
$$
\begin{aligned}
&\int^{T}_{0}\g|\va\n A(\rho)|^{2}\Lambda^{s}(\va C(\rho))\leq {\cal E}_{0}^{\frac{1}{2}}
\|\Lambda^{1+\frac{s}{2}}(\va B(\rho))\|_{L^{2}(L^{2})}^{2}.
\end{aligned}
$$
For the other term 
we proceed similarly as 
theorem \ref{3T1}.
\hfill {$\Box$}
\subsubsection{What happens when we choose a coefficient of capillarity $\kappa$ large}
Assume that $\kappa(\rho)=\kappa\rho^{\alpha}$ with $\alpha>0$ and $\kappa$ is a very large constant coefficient. We have then by energy inequality:
$$\|\n A(\rho)\|_{L^{\infty}(L^{2})}\leq C(\|\n A(\rho_{0})\|_{L^{2}}+\frac{1}{\kappa}(\|\sqrt{\rho_{0}}u_{0}\|_{L^{2}}+\|\rho\|_{L^{\gamma}}^{\frac{\gamma}{2}})\big).$$
It means that we can hope gain of derivative on the density if we assume that $\kappa$ is enough big and only $\n A(\rho_{0})$
small in $L^{2}$. We obtain then the following theorem:
\begin{theorem}
Let $N=2$ and $\alpha,M\in\R$. $(\rho,u)$ is a smooth approximate solution of the system
$(\ref{3systeme})$ with initial data in the energy space and with the following capillarity coefficient:
$$\kappa(\rho)=\frac{1}{\rho^{2}}1_{\{\rho<\alpha\}}+\theta_{1}(\rho)1_{\{\alpha\leq\rho\leq2
\alpha\}}+\kappa 1_{\{2\alpha<\rho\}}.$$
where $\theta_{1},\;\theta_{2}$ are regular function such that $\kappa$ is a regular function.
Then there exists a constant
$\eta>0$ and $\kappa$ enough big depending only on the constant coming from the Sobolev
embedding and on the initial data such that if:
$$\|\n A(\rho_{0})\|_{L^{2}(\R^{2})}\leq\eta$$
then we get for all $\va\in C^{\infty}_{0}(\R^{N})$:
$$\|\va B(\rho)\|_{L^{2}_{T}(\dot{H}^{1+\frac{s}{2}})}\leq M\;\;\;\;\mbox{with}\;\;0<s<2,$$
where $M$ depends only on the initial conditions data, on
$T$, on $\va$ and on $s$.
\end{theorem}
{\bf Proof:} In fact by following the proof of theorem \ref{3T1}, we just have to check that:
$$\kappa\int^{T}_{0}\int_{\R^{N}}\va|\n A(\rho)|^{2}\Lambda^{s}(\va B(\rho))dxdt<\kappa\|B(\rho)\|_{L^{2}(\dot{H}^{1+\frac{s}{2}})}^{2}.$$
And this is the case because with our hypothesis:
$$\|\n A(\rho)\|_{L^{\infty}(L^{2})}\leq C(\|\n A(\rho_{0})\|_{L^{2}}+\frac{1}{\kappa}(\|\sqrt{\rho_{0}}u_{0}\|_{L^{2}}+\|\rho\|_{L^{\gamma}}^{\frac{\gamma}{2}})\big)\leq\e,$$
with $\e$ enough small when $\kappa$ is enough large and $\|\n A(\rho_{0})\|_{L^{2}}$ enough small.
\hfill {$\Box$}
\subsection{Control of the density without any conditions of smallness on the initial data}
In this section, we want avoid to impose some conditions of smallness on the initial data. We would get a theorem of global weak solution with large initial data. For making, we have to localize our previous arguments with test function $\va\in C^{\infty}_{0}$
with small diameter support. Indeed we would get a control of $\|1_{B(x,r)}\n A(\rho)\|_{L_{T}^{\infty}(L^{2})}$ and prove that this norme is small uniformly for all $x\in \R^{N}$  when $r$ is small.
So for a chosen compact $K$, we have to split it in a finite union of small ball and apply our previous argument on each of these balls.
Before entering in the heart of subject, we would start now with localizing the classical energy inequality on small support.
\subsubsection{Localizing energy inequality and phenomena of concentration}
We are interesting in getting new energy inequalities which take in account the support of the initial data. In particular we want investigate what happend when we choose initial data localize in small ball. So for making we let $\va\in C^{\infty}_{0}(\R^{N})$ with $\va=1$ on $B(x_{0},R)$ and $\mbox{supp}\va\subset B(x_{0},2R)$. Multiplying the equation of momentum conservation in the system
(\ref{3systeme}) by $\va u$ and integrating by parts over $(0,t)\times\R^{N}$,
we obtain the following
estimate:
\begin{equation}
\begin{aligned}
&\frac{1}{2}\int_{\R^{N}}\va(x)\big(\rho|u|^{2}(t,x)+\frac{\kappa(\rho)}{2}|\n\rho|^{2}+2\Pi(\rho)-2\Pi(\bar{\rho})\big)(t,x)dx\\
&+\int_{0}^{t}\int_{\R^{N}}(\mu|\n u|^{2}+\xi|{\rm div}u|^{2})\va dx\leq
\int_{\R^{N}}\frac{1}{2}\va(x)\big(\rho_{0}u_{0}^{2}+\kappa(\rho_{0})|\nabla\rho_{0}|^{2})dx\\
&+\g\va(x)(\Pi(\rho_{0})-\Pi(\bar{\rho}))(x)dx+\int_{0}^{t}\int_{\R^{N}}\big[(u\cdot\n u).\n\va +P(\rho)u\cdot\n\va\big]dxdt\\
&+\int_{\R}\g\kappa(\rho)\n\rho\cdot\n\va{\rm div}(\rho u) \,dx\,dt+
\int_{\R}\g\big(\kappa(\rho)+\frac{1}{2}\rho\kappa^{'}(\rho)\big)|\n\rho|^{2}u\cdot\n\va dxdt\\
&+\int_{\R}\g\rho\kappa(\rho)\p_{j}\rho\p_{j}u_{i}\p_{i}\va\,dxdt+\int_{\R}\g\rho\kappa(\rho)\p_{j}\rho u_{i}\p_{i,j}\va,dxdt.
\end{aligned}
\label{3localize}
\end{equation}
Now we have to control the term on the right handside in the goal to get some energy estimate localized in space. We recall that in the proof of theorem \ref{3T1}, we need a condition of smallness on $\va\n A(\rho)\in L^{\infty}(L^{2})$. So by a condition of small support on $\va$ we can hope get this condition of smallness on $\va\n A(\rho)\in L^{\infty}(L^{2})$ for a finite time. In fact we search to prove that there is not concentration effect in
some point of the space for a small intervall $[0,T]$. More precisely it means that for any sequel of approximate solution $(\rho_{n},u_{n})_{n\in\mathbb{N}}$, we have the following property:
$$\forall K\;\;\mbox{a compact}, \exists\e>0\;\forall t\in[0,T]\Rightarrow
\|\n \big(A(\rho_{n})(t,\cdot)\big)1_{K}(\cdot)\|_{L^{2}}\leq\e.$$
So we need of this type of property to can use a boobstrap in the proof of theorem \ref{3T1}
without asking condition of smallness on the initial data. Moreover this property is very natural. Indeed we want avoid some concentration effect. In particular if $(\rho_{n},u_{n})$ is a sequel of approximate solution for the system (\ref{3systeme}), we want ask that $|\n A(\rho_{n})|^{2}$ doesn'tconverge to some Dirac measures.
\begin{proposition}
Let $\alpha>0$, $\beta>0$,$\eta>0$, $\kappa>0$, $\theta$ a regular function such that:
$$\kappa(\rho)=\frac{1}{\rho^{2}}1_{\{\rho<\alpha\}}+\theta (\rho)1_{\{\alpha\leq\rho<2\alpha\}}+\kappa 1_{\{2\alpha\leq\rho\leq\eta\}}+\frac{1}{\rho^{2}}1_{\{\rho\geq\eta\}}.$$
Let $(\rho,u)$ a regular approximate solution of system (\ref{3systeme}) with large initial data in the energy space.
Assuming that $u\in L^{1+\beta}_{T}(L^{\infty})\cap L^{2+\beta}(L^{2})$ then it exist a time $T_{0}>0$ such that for all $\va\in C^{\infty}_{0}(\R^{N})$ we have:
$$\|\va B(\rho)\|_{L^{2}_{T_{0}}(H^{1+\frac{s}{2}})}\leq M,$$
where $M$ depends only of the initial data, $K$ and $\| u\|_{L^{1+\beta}_{T}(L^{\infty})}$.
\end{proposition}
{\bf Proof:}
The proof consists only to get a control of $\|\n \big(A(\rho)\big)1_{K}\|_{L^{\infty}(L^{2})}$
for $K$ a compact of small enough measure. Indeed by following the same lines
as the proof of theorem \ref{3T1}, the main difficulty is to control the term for $\va\in C^{\infty}_{0}(\R^{N})$:
$$\int_{0}^{T}\g\va|\n A(\rho)|^{2}\Lambda^{s}B(\rho)dxdt.$$
In our case, we have to take in consideration that the integration is only on the support of $\va$.
So let $K=\mbox{supp}\va$ a compact of $\R^{N}$, we 
consider a subordinate partition of unity of the indicatrice $1_{K}$ with $(\va_{k})_{1\leq k\leq p}$ of class $C^{\infty}$ such that:
\begin{itemize}
\item $\mbox{supp}\va_{k}\subset B(x_{k},\lambda)\subset K$,
    \;\;$\sum_{k}\va_{k}=1$ on $K$ and $0\leq\va_{k}\leq 1$,
\item $\|\n^{\alpha}\va_{k}\|_{L^{\infty}}\leq C_{\alpha}|\lambda|^{-\alpha}$ for $|\alpha|\leq2$,
\item $\|\n\va_{k}\|_{L^{p}}\leq C_{\alpha}\frac{1}{|\lambda|^{1-\frac{1}{p}}}$ for $1\leq p\leq\infty$.
\end{itemize}
We have now to estimate the following term:
\begin{equation}
\begin{aligned}
&\int_{0}^{T}\g\va_{k}|\n A(\rho)|^{2}\Lambda^{s}(\va_{k}B(\rho))dxdt=\\
&\int_{0}^{T}\g\n A(\rho)\cdot(\va_{k}\frac{A^{'}(\rho)}{B^{'}(\rho)}\n B(\rho))\Lambda^{s}(\va_{k}B(\rho))dxdt+R_{\va_{k}}
\end{aligned}
\label{3crucialpoint}
\end{equation}
to control $\va_{k}B(\rho)$ in $L^{2}(H^{1+\frac{s}{2}})$ and where here $R_{\va_{k}}$ is a terme easy to control. We want use a argument of boobstrap as in the proof of theorem \ref{3T1}.
We have now just to prove that  for $\e$ enough small, it exists a $T_{0}$ such that we have:
$$\|1_{\mbox{supp}\va_{k}}\n A(\rho)\|_{L^{\infty}_{T_{0}}(L^{2})}\leq\e.$$
Let $\psi_{k}\in C^{\infty}_{0}$ such that $\psi_{k}=1$ on $\mbox{supp}\va_{k}$ and $\mbox{supp}\psi_{k}\subset 2\mbox{supp}\va_{k}$
We set then:
$$A(T,\psi_{k})=\frac{1}{2}\int_{\R^{N}}\psi_{k}(x)\big(\rho|u|^{2}(t,x)+|\n A(\rho)|^{2}+2\Pi(\rho)-
2\Pi(\bar{\rho})\big)(t,x)dx.$$
By the equation (\ref{3localize}), we have:
$$
\begin{aligned}
&A(T,\psi_{k})+\int_{0}^{t}\int_{\R^{N}}(\mu|\n u|^{2}+\xi|{\rm div}u|^{2})\psi_{k} dx\leq
\int_{\R^{N}}\frac{1}{2}\psi_{k}(x)\big(\rho_{0}u_{0}^{2}+|\nabla A(\rho_{0})|^{2})dx\\
&+\g\psi_{k}(x)(\Pi(\rho_{0})-\Pi(\bar{\rho}))(x)dx+\int_{0}^{t}\int_{\R^{N}}\big[(u\cdot\n u).\n\psi_{k} +P(\rho)u\cdot\n\psi_{k}\big]dxdt\\
&+\int_{\R}\g\kappa(\rho)\n\rho\cdot\n\psi_{k}{\rm div}(\rho u) \,dx\,dt+
\int_{\R}\g\big(\kappa(\rho)+\frac{1}{2}\rho\kappa^{'}(\rho)\big)|\n\rho|^{2}u\cdot\n\psi_{k} dxdt\\
&+\int_{\R}\g\rho\kappa(\rho)\p_{j}\rho\p_{j}u_{i}\p_{i}\psi_{k}\,dxdt+\int_{\R}\g\rho\kappa(\rho)\p_{j}\rho u_{i}\p_{i,j}\psi_{k},dxdt.
\end{aligned}
$$
We have to control all the term of right to have estimates on $A(T,\psi_{k})$. Let $\e>0$, easily for $\lambda$ enough big, we have:
\begin{equation}
\int_{\R^{N}}\frac{1}{2}\psi_{k}(x)\big(\rho_{0}u_{0}^{2}+|\nabla A(\rho_{0})|^{2})dx
+\g\psi_{k}(x)(\Pi(\rho_{0})-\Pi(\bar{\rho}))(x)dx<\frac{\e}{2}.
\label{initialdata}
\end{equation}
We have now to control the other terms, we begin by the most complicated and we have:
$$
\begin{aligned}
&\int_{\R}\g\kappa(\rho)\n\rho\cdot\n\psi_{k}{\rm div}(\rho u) \,dx\,dt\leq\|\n A(\rho)\|_{L^{\infty}(L^{2})}^{2}\|\n\psi_{k}\|_{L^{\infty}}\|u\|_{L^{1+\beta}(L^{\infty})}T^{\frac{\beta}{1+\beta}}\\
&\hspace{7cm}+\|\n A(\rho)\|_{L^{\infty}(L^{2})}\|\n\psi_{k}\|_{L^{\infty}}\|\n u\|_{L^{2}(L^{2})}T^{\frac{1}{2}}.
\end{aligned}
$$
We can treat similarly the other terms and we get finally by using inequality (\ref{initialdata}):
$$
\begin{aligned}
&A(T,\psi_{k})\leq\frac{\e}{2}+\|\n A(\rho)\|_{L^{\infty}(L^{2})}^{2}\|\n\psi_{k}\|_{L^{\infty}}\|u\|_{L^{1+\beta}(L^{\infty})}T^{\frac{\beta}{1+\beta}}\\
&+
\|A(\rho)\|_{L^{\infty}(\dot{H}^{1})}\|\n\psi_{k}\|_{L^{\infty}}\|\n u\|_{L^{2}(L^{2})}T^{\frac{1}{2}}+\| A(\rho)\|_{L^{\infty}(\dot{H}^{1})}\|D^{2} \psi_{k}\|_{L^{2-\beta}}\|\n u\|_{L^{2}(L^{2})}T^{\frac{1}{2}}\\
&+\|u\|_{L^{2}(\dot{H}^{1})}\|u\|_{L^{2+\beta}(L^{2})}\|\n\psi_{k}\|_{L^{\infty}}T^{\frac{\beta}{2(2+\beta)}}+
\|P(\rho)\|_{L^{\infty}(L^{1})}\|u\|_{L^{1+\beta}(L^{\infty})}\|\n\psi_{k}\|_{L^{\infty}}T^{\frac{\beta}{1+\beta}},
\end{aligned}
$$
By the conditions on $\va_{k}$ and the previous inequality we get:
$$
\begin{aligned}
&A(T,\psi_{k})\leq \frac{\e}{2}+C|\lambda|\big(\|u\|_{L^{1+\beta}(L^{\infty})}T^{\frac{\beta}{1+\beta}}+T^{\frac{1}{2}}+
\|u\|_{L^{2+\beta}(L^{2})}T^{\frac{\beta}{2(2+\beta)}}\\
&\hspace{8cm}+\|u\|_{L^{1+\beta}(L^{\infty})}T^{\frac{\beta}{1+\beta}}+|\lambda|T^{\frac{1}{2}}\big).
\end{aligned}
$$
For a small enough time $T_{0}$ depending of the initial data, $\|u\|_{L^{1+\beta}(L^{\infty})}$, $\|u\|_{L^{2+\beta}(L^{2})}$ and
$\lambda$ we have:
\begin{equation}
A(T,\psi_{k})\leq\e.
\label{3important}
\end{equation}
We can now come back to the crucial point of the proof of theorem \ref{3T1}, and from the equality (\ref{3crucialpoint})
and inequality (\ref{3important}) we have:
$$
\begin{aligned}
\int_{0}^{T}\g\va_{k}|\n A(\rho)|^{2}\Lambda^{s}(\va_{k}B(\rho))dxdt&\leq C\|\psi_{k}\n A(\rho)\|_{L^{\infty}(L^{2})}\|\va_{k}B(\rho)\|^{2}_{L^{2}(H^{1+\frac{s}{2}})}+M_{0},\\
&\leq C\e\|\psi_{k}\n A(\rho)\|_{L^{\infty}(L^{2})}\|\va_{k}B(\rho)\|^{2}_{L^{2}(H^{1+\frac{s}{2}})}+M_{0}.
\end{aligned}
$$
By choosing $\e$ enough small, we conclude that $\va_{k}B(\rho)\in L^{2}(H^{1+\frac{s}{2}})$. Next we get by Sobolev embedding
$\va_{k}\n B(\rho)\in L^{2}_{T_{0}}(L^{p})$ with $\frac{1}{p}=\frac{1}{2}-\frac{s}{4}$. So we have $1_{K}\n B(\rho)=\sum_{k}\va_{k}\n B(\rho)$ which is in $L^{2}(L^{p})$, by coming back at the proof of theorem \ref{3T1} we finally get a control of
$\va B(\rho)$ in $L^{2}(H^{1+\frac{s}{2}}$ and this conclude the proof.
\hfill {$\Box$}
\begin{remarka}
In a similar way, we could show that if we control the high frequencies of $\n A(\rho)$ in $L^{\infty}(L^{2})$
then we can get a gain of derivative with large initial data in the energy space. More precisely we have:
$$\hat{\p_{i}A(\rho)}=\hat{\p_{i}A(\rho)}1_{\{|\xi|\leq M\}}+\hat{\p_{i}A(\rho)}1_{\{|\xi|>M\}}=f_{1}+f_{2}.$$
We know that $f_{1}$ is regular and in particular $f_{1}\in L^{\infty}(L^{p})$ with $p\in[2,+\infty[$, so by using the same argument
as in proff \ref{3T1}, we are able to treat this term.
The main difficulty is to control $\hat{\p_{i}A(\rho)}1_{\{|\xi|>M\}}$, indeed we would use a argument of smallness on
$\|\hat{\p_{i}A(\rho)}1_{\{|\xi|>M\}}\|_{L^{\infty}(L^{2})}$ for $M$ enough big.\\
In particular for an approximate sequence of solutions, is it possible to find for a $\e>0$ $M$ enough big and depending only on the initial data
such that:
$$
\begin{aligned}
&\|\hat{\p_{i}A(\rho_{n})}1_{\{|\xi|>M\}}\|_{L^{\infty}(L^{2})}\leq\e.\\
\end{aligned}
$$
The difficulty is to show that the mass does not run away to the inifinity.
\end{remarka}
\subsection{Existence of global weak solutions for $N=2$}
We may now turn to our compactness result. First, we assume that a
sequence $(\rho_{n},u_{n})_{n\in\mathbb{N}}$ of approximate weak
solutions has been constructed by a mollifying process, which have
suitable regularity to justify the formal
estimates like the classical energy estimates and these coming from your previous theorems. 
\begin{remarka}
In fact,
we easily construct this sequence of solution $(\rho_{n},u_{n})_{n\in\mathbb{N}}$ by using the result
of \cite{3H}, it suffices to choose initial data $(\rho^{0}_{n},u^{0}_{n})\in B^{\N}_{2,1}\times B^{\N-1}_{2,1}$
with a condition of smallnees. We know indeed that there exists $\eta>0$ such that if:
$$\|\rho_{n}^{0}\|_{B^{\N}_{2,1}}+\|u_{n}^{0}\|_{B^{\N_1}_{2,1}}\leq\eta,$$
then the solution $(\rho^{0}_{n},u^{0}_{n})$ are global and strong.
\label{remarqueexistence}
\end{remarka}
Moreover in the sequel the viscosity coefficients check the properties of definition \ref{3defviscosity}, it means that they are non degenerate. Moreover in the sequel we will work in dimension $2$, and to simplify we choose a isentropic pressure $P(\rho)=\rho^{\gamma}$ with $\gamma>1$.\\
Moreover this sequence  $(\rho_{n},u_{n})_{n\in\mathbb{N}}$ has
initial data $((\rho_{0})_{n},(u_{0})_{n}))$ close to the energy
space. By using the above energy inequalities, we assume that
$j_{\gamma}((\rho_{0})_{n})$, $|\n(A\big((\rho_{0})_{n}\big)|^{2}$ and
$(\rho_{0})_{n}|(u_{0})_{n}|^{2}$ are bounded in
$L^{1}(\R^{N})$ so that $(\rho_{0})_{n}$ is bounded in $L^{\gamma}_{2}(\R^{N})$ .\\
\\
Then it follows from the energy inequality that:
\begin{enumerate}
\item $j_{\gamma}(\rho_{n}),|\n A(\rho_{n})|^{2},\,\rho_{n}|u_{n}|^{2}$ are bounded uniformly in $L^{\infty}(0,T,L^{1}(\R^{N}))$,
\item $Du_{n}$ is bounded uniformly in $L^{2}(\R^{N}\times
(0,T))$,
\item $u_{n}$ is bounded uniformly
in $L^{2}(0,T,H^{1}(B_{R}))$ for all $R,T\in (0,+\infty)$.
\end{enumerate}
Extracting subsequences if necessary, we may assume that
$\rho_{n},\,u_{n}$ converge weakly respectively in
$L^{\infty}(L^{\gamma}_{2}(\R^{N}))$, $L^{2}(0,T;H^{1}(B_{R}))$ to
$\rho,\,u$ for all $R,T\in(0,+\infty)$. In a similar way to the case of Navier-Stokes compressible, we can extract subsequences
such that
$\sqrt{\rho_{n}}u_{n},\,\rho_{n}u_{n},\,\rho_{n}u_{n}\otimes u_{n}$
converge weakly to $\sqrt{\rho}u$, $\rho u$ and $\rho u\otimes u$.\\
In fact the main difficulty is to verify that the quadratic gradient term of the density $\n A(\rho_{n})\otimes\n A(\rho_{n})$
converge to $\n A(\rho)\otimes\n A(\rho)$.\\
Finally we add uniform hypothesis on $B(\rho_{n})$ coming from the gain of regularity that we have obtained in the previous part, so we assume that:
$$
\begin{aligned}
&\forall\va\in C_{0}^{\infty}(\R^{N}),\;\;\;\;\va B(\rho_{n})\;\;
\mbox{is bounded uniformly in}\;\; L^{2}_{T}(H^{1+\frac{s}{2}})\:\:\:\:\mbox{with}\;\;s>0.
\end{aligned}
$$
Moreover when $\kappa(\rho)=\frac{\kappa}{\rho^{2}}$ with $\kappa>0$ we have:
\begin{itemize}
\item $\va\rho_{n}$ is bounded uniformly in $L^{\gamma+\alpha}
((0,T)\times\R^{N})$.
\item $\va\rho_{n}^{\alpha-2}\n\rho_{n}$
is bounded uniformly in $L^{2}_{T}(L^{2}(\R^{N}))$.
\end{itemize}
We can now show the two following theorem. The next theorem concerns the existence of global weak solutions in the case where $\kappa(\rho)=\kappa\rho^{-2}$.
\begin{theorem}
\label{3Tfinal1}
Here we assume 
$\kappa(\rho)=\frac{\kappa}{\rho^{2}}$ with $\kappa>0$.
There exists $\eta>0$ such that if:
$$\|\n(\ln\rho_{0}^{n})\|_{L^{2}}+\|\sqrt{\rho_{0}^{n}}|u_{0}^{n}|\|_{L^{2}}
+\|j_{\gamma}(\rho_{0}^{n})\|_{L^{1}}\leq\eta$$ then, up to an extraction a
subsequence $(\rho_{n},u_{n})$ converges strongly to a weak solution
 $(\rho,u)$ (see definition \ref{3defexistence}) of the system (\ref{3systeme}). Moreover we
have $\n \log(\rho_{n})\otimes\n \log(\rho_{n})$ converges strongly in
$L^{1}_{loc}(\R\times\R^{N})$. In addition $\rho$ check for all $\va\in C^{0}_{\infty}$:
$$\|\va \ln(\rho) \|^{2}_{L^{2}_{T}(H^{1+\frac{s}{2}})}+\|\va\rho^{\alpha-2}\n\rho\|^{2}_{L^{2}(L^{2})}+\|\rho\|_{L^{\gamma+\alpha}
((0,T)\times\R^{N})}^{\gamma+\alpha}\leq M\;\;\;\;\;\mbox{with}\;\;0\leq
s\leq 2,$$
where $M$ depends only on the initial conditions data, on $T$, on $\va$ and on $s$.
\end{theorem}
\begin{remarka}
This theorem is a theorem of existence of global weak solution and not only a result of stability.
Indeed as explained in the remarka \ref{remarqueexistence} we are able by using the result of \cite{3H}
to construct global approximate solution. And the result of \cite{3H} is compatible with our hypothesis of smallness.
It would be the same for the following theorems.
\end{remarka}
The next theorem treat of global weak solution for capillarity coefficients approximating a constant.
\begin{theorem}
\label{3Tfinal2} Let $N=2$, $\alpha>0$, $\e\geq0$ and the following capillary coefficient:
$$\kappa(\rho)=\frac{1}{\rho^{2+\e}}1_{\{\rho<\alpha\}}+\theta_{1}(\rho)1_{\{\alpha\leq\rho\leq2
\alpha\}}+\kappa 1_{\{\rho>2\alpha\}}.$$
where $\theta_{1},\;\theta_{2}$ are regular function ssuch that $\kappa$ is a regular function.
There exists $\eta>0$ such that if:
$$\|\n(A(\rho_{0}^{n}))\|_{L^{2}}+\|\sqrt{\rho_{0}^{n}}|u_{0}^{n}|\|_{L^{2}}
+\|j_{\gamma}(\rho_{0}^{n})\|_{L^{1}}\leq\eta$$ then, up to a
subsequence $(\rho_{n},u_{n})$ converges strongly to a weak solution
 $(\rho,u)$ (see definition \ref{3defexistence}) of the system (\ref{3systeme}). Moreover we
have $\n A(\rho_{n})\otimes\n A(\rho_{n})$ converges strongly in
$L^{1}_{loc}(\R\times\R^{N})$.
In addition $\rho$ check for all $\va\in C^{0}_{\infty}$:
$$\|\va A(\rho) \|^{2}_{L^{2}_{T}(H^{1+\frac{s}{2}})}\leq M\;\;\;\;\;\mbox{with}\;\;0\leq
s\leq 2,$$
where $M$ depends only on the initial conditions data, on $T$, on $\va$ and on $s$.
\end{theorem}
To finish, we give a theorem of existence of global weak solution when the capillarity coefficient $\kappa$ is big and when only $\n A(\rho_{0})$ admits a condition of smallness.
\begin{theorem}
\label{3Tfinal3} Let $N=2$, $\alpha>0$, $\e\geq0$ and the following capillary coefficient:
$$\kappa(\rho)=\frac{1}{\rho^{2+\e}}1_{\{\rho<\alpha\}}+\theta_{1}(\rho)1_{\{\alpha\leq\rho\leq2
\alpha\}}+\kappa 1_{\{\rho>2\alpha\}}.$$
where $\theta_{1},\;\theta_{2}$ are regular function ssuch that $\kappa$ is a regular function.
Let $\sqrt{\rho_{0}}u_{o}\in L^{2}$, $j_{\gamma}(\rho_{0}\in L^{1}$
There exists $\kappa$ enough big and $\eta>0$ such that if:
$$\|\n(A(\rho_{0}^{n}))\|_{L^{2}}\leq\eta$$ then, up to a
subsequence $(\rho_{n},u_{n})$ converges strongly to a weak solution
 $(\rho,u)$ (see definition \ref{3defexistence}) of the system (\ref{3systeme}). Moreover we
have $\n A(\rho_{n})\otimes\n A(\rho_{n})$ converges strongly in
$L^{1}_{loc}(\R\times\R^{N})$.
In addition $\rho$ check for all $\va\in C^{0}_{\infty}$:
$$\|\va A(\rho) \|^{2}_{L^{2}_{T}(H^{1+\frac{s}{2}})}\leq M\;\;\;\;\;\mbox{with}\;\;0\leq
s\leq 2,$$
where $M$ depends only on the initial conditions data, on $T$, on $\va$ and on $s$.
\end{theorem}
{\bf Proof of the theorem \ref{3Tfinal1} :}
The main difficulty states in proving the convergence of the following nonlinear terms $P(\rho_{n})$ and
$\n\ln\rho_{n}\otimes\n\ln\rho_{n}$. The other terms follows the same lines as the compressible Navier-Stokes problem
studied by P-L. Lions, E. Feireisl, A. Novot\'y and many other authors.\\
According to theorem \ref{3T1} we have seen that for all $\va\in
C^{\infty}_{0}(\R^{N})$ $\va\ln\rho_{n}\in
L^{2}_{T}(H^{1+\frac{s}{2}})$.
So we have $\forall\va\in C^{\infty}_{0}(\R^{N})$
\\
We can now use some results of compactness to show that
$\n(\ln\rho_{n})$ converge strongly in $L^{2}_{T}(L^{2}_{loc})$ to
$\n\ln\rho$. We recall the following theorem from Aubin-Lions ( see
Simon for general results \cite{3Sim}).
\begin{lemme}
\label{3Aubin} Let $X\hookrightarrow Y\hookrightarrow Z$ be Hilbert
spaces such that the embedding from $X$ in $Y$ is compact. Let
$(f_{n})_{n\in\mathbb{N}}$ a sequence bounded in $L^{q}(0,T;Y)$,
(with $1<q<+\infty$) and $(\frac{d f_{n}}{dt})_{n\in\mathbb{N}}$
bounded in $L^{p}(0,T;Z)$ (with $1<p<+\infty$), then
$(f_{n})_{n\in\mathbb{N}}$ is relatively compact in $L^{q}(0,T;Y)$.
\end{lemme}
We need now to localize our arguments because we want use some result of compactness for
the local Sobolev space $H_{loc}^{\frac{s}{2}}$ with $s>0$. $H_{loc}^{\frac{s}{2}}$ is compactly
embedded in $L^{2}_{loc}$. Let $(\chi_{p})_{p\in\mathbb{N}}$ be a
sequence of $C^{\infty}_{0}(\R^{N})$ cut-off functions supported in
the ball $B(0,p+1)$ of $\R^{N}$ and equal to 1 in a neighborhood of
$B(0,p)$.
We have then by using mass equation:
$$\frac{d}{dt}\n(\ln\rho_{n})+\n{\rm
div}(u_{n})+\n(u_{n}\cdot\n\ln\rho_{n})=0$$ We can
then show that
$\big(\frac{d}{dt}(\chi_{p}\n(\ln\rho_{n}))\big)_{n\in\mathbb{N}}$
is uniformly bounded for all $p$ in $L^{q}_{T}(H^{\alpha})$ for
$\alpha<0$ by using energy inequalities. Moreover
$\big(\chi_{p}\n(\ln\rho_{n})\big)_{n\in\mathbb{N}}$ is uniformly
bounded for all $p$ in $L^{2}_{T}(H^{\frac{s}{2}})$. Applying lemma
\ref{3Aubin} with the family
$\big(\chi_{p}\n(\ln\rho_{n})\big)_{n\in\mathbb{N}}$ and
$X=\chi_{p}H^{\frac{s}{2}}$, $Y=\chi_{p}L^{2}$, $Z=\chi_{p}
H^{\alpha}$ and using Cantor's diagonal process, we
provides that after up to a subsequence:
\begin{equation}
\forall
p>0\;\;\;\chi_{p}\n(\ln\rho_{n})\rightarrow_{n\rightarrow+\infty}
\chi_{p}a\;\; \mbox{in}\;\;L^{2}_{T}(L^{2}). \label{3egalite}
\end{equation}
with for all $p\in\mathbb{N}$ , $\chi_{p}a\in L^{2}(H^{\frac{s}{2}})$.\\
Moreover as
$\frac{1}{2}\sqrt{\rho_{n}}\n\rho_{n}=\n\sqrt{\rho_{n}}\in L^{\infty}_{T}(L^{p})$ with $\frac{1}{p}=\frac{1}{2}+\frac{1}{2\gamma}$, by the same argument we have $\sqrt{\rho_{n}}$ converge strongly
to a certain $b$ in $L^{2}(L_{loc}^{2\gamma-\e})$ with $\e>0$ such that $2\gamma-\e=2$. Moreover up a subsequence
$\sqrt{\rho_{n}}$ converges a.e to $b$.
We have then for all
$\va\in C^{\infty}_{0}$:
$$\int_{(0,T)\times\R^{N}}\va\rho^{n}dxdt\rightarrow_{n\rightarrow+\infty}\int_{(0,T)\times\R^{N}}\va b^{2}dxdt.$$
And as $\rho_{n}$ converges weakly to $\rho$, we have shown that $b^{2}=\rho$ and that $\rho_{n}$ converges a.e to
$\rho$.\\
We can show that for all $\va\in C^{\infty}_{0}$ $\va\ln\rho_{n}$ converges weakly to $\va\ln\rho$ by the fact that
$\rho_{n}$ converges a.e to
$\rho$ and the proposition \ref{3propzuily}.
It means that $a=\n\ln\rho$.
Finally we have shown that $\n(\ln\rho_{n})$ converges strongly to $\n\ln\rho$ in $L^{2}_{T}(L^{2}_{loc})$.
We have then obtained that $\n\ln\rho_{n}\otimes\n\ln\rho_{n}$ converges in distribution sense to
$\n\ln\rho\otimes\n\ln\rho$.
\\
The last difficulty is to treat the term $P(\rho^{n})$, we proceed similarly. Let $\va\in C^{\infty}_{0}(\R^{N})$, so
as $\va \n \sqrt{\rho_{n}}=\frac{1}{2}\va\sqrt{\rho_{n}}\n\ln\rho_{n}$ and as $\rho_{n}$ is uniformly bounded in $ L^{\gamma+\alpha}_{T}(L^{\gamma+\alpha}_{loc})$ we get $\n\sqrt{\rho_{n}}$  is uniformly bounded in $L^{2(\gamma+\alpha)}_{T}(L_{loc}^{p})$ with $\frac{1}{p}=\frac{1}{2}+\frac{1}{2(\gamma+\alpha)}$. We can conclude by compact Sobolev embedding and proposition \ref{3propzuily}.
\hfill {$\Box$}\\
The proof of theorem \ref{3Tfinal2} and \ref{3Tfinal3} follows the same line than the previous proof.
\section{Existence of weak solution in the case $N=1$}
\label{3S3} We are now interested by the case $N=1$. To start with,
we focus on the gain of derivative for $\n A(\rho)$ with a general capillarity term.
\subsection{Gain of derivative}
We can now write a theorem where we reach a gain of derivative on
the density $\rho$ by using the same type of inequalities as in the
case $N=2$ and without any conditions of smallness on the initial data.
\begin{theorem}
Let $(\rho,u)$ be a regular solution of the system (\ref{3systeme})
with initial data in the energy space and let $\kappa(\rho)=\kappa\rho^{\alpha}$ a general capillarity coefficient with $\alpha\in\R$. Then we have:
$$\|B(\rho)\|_{L^{2}_{T}(H^{1+\frac{s}{2}}(\R))}\leq M_{0}$$
with $0\leq s<\frac{1}{2}$ and $M_{0}$ depending only of the initial
data. \label{3T3}
\end{theorem}
\begin{remarka}
We observe the two important facts:
\begin{enumerate}
\item We don't need any hypothesis on the size of the initial
data.
\item We don't need to localize because we know that $\rho\in L^{\infty}_{t,x}$.
\item We don't need to assume that $\frac{1}{\rho}\in L^{\infty}$.
\end{enumerate}
\end{remarka}
{\bf Proof of theorem \ref{3T3} :}\\
\\
We use here the same estimates as in the previous proof except for the
delicate term:
$\int_{0}^{T}\int_{\R}|\p_{x} A(\rho)|^{2}\La^{s}B(\rho)$.
In the sequel we will show only the case $\alpha\geq -2$, the proof of the other case is similar. We have then $\p_{x}A(\rho)\in L^{\infty}(L^{2})$ and
$\rho-\bar{\rho}\in L^{\infty}(L^{2})$ so 
by Sobolev embedding $A(\rho)\in
L^{\infty}(L^{\infty})$. 
Next we set $\n B(\rho)=\frac{B^{'}(\rho)}{A^{'}(\rho)}\n A(\rho)=\frac{1}{\kappa}A(\rho)\n A(\rho)$.
But we know that $A(\rho)$ belongs to $L^{\infty}(L^{\infty})$, so that $\n B(\rho)\in L^{\infty}(L^{2})$.
Finally we get for $0<s\leq1$,
$\La^{s}B(\rho)\in L^{\infty}(H^{1-s})$. Now for $0\leq
s<\frac{1}{2}$ by Sobolev embedding we obtain:
$$\La^{s}B(\rho)\in L^{\infty}(L^{\infty}).$$
So we can control the term
$\int_{0}^{T}\int_{\R}|\p_{x}A(\rho)|^{2}|\La^{s}B(\rho)$ as follows:
$$\int_{0}^{T}\int_{\R}|\p_{x}A(\rho)|^{2}|\La^{s}B(\rho)|\lesssim\|\p_{x}A(\rho)\|_{L^{\infty}_{T}(L^{2})}^{4}.$$
We treat the other terms similarly as in the previous proof. \hfill{$\Box$}
\subsection{Results of compactness}
We can now prove our result of stability of solution in the case
$N=1$ by using the previous gain of derivative. Let
$(\rho_{n},u_{n})_{n\in\mathbb{N}}$ a sequel of approximate weak
solutions of system (\ref{3systeme}). We get now a theorem of global weak solution with large initial data.
\begin{theorem}
Let $(\rho_{0}^{n},u_{0}^{n})$ initial data of the system
(\ref{3systeme}) in the energy space 
Let $\kappa(\rho)=\frac{1}{\rho^{2+\e}}1_{\{\rho\leq\alpha\}}+\theta(\rho)1_{\{\alpha<\rho\leq 2\alpha\}}+\kappa
1_{\{<\rho>2\alpha\}}$ with $\e\geq 0$ and $\theta$ a regular function.\\ Then 
up to a subsequence, $(\rho_{n},u_{n})$ converges
strongly to a weak solution $(\rho,u)$
on $(0,T)\times\R$ in the sense of the distribution  for all $T\in(0,+\infty)$ (see definition \ref{3defexistence}).
Moreover $\p_{x}A(\rho_{n})$ converges strongly in
$L^{2}(0,T,L^{2}_{loc}(\R^{N}))$ to $\p_{x}A(\rho)$. \label{3T4}
\end{theorem}
{\bf Proof:}The proof follows the same lines as in the proof of theorem \ref{3Tfinal1}.
\hfill{$\Box$}
\section{Appendix}
\label{appendix}
\subsection{Computation of the capillary term}
This section is devoted to rewrite clearly the capillarity tensor $K$ in the the goal to
express the non linear terms in distribution
sense. On other interesting reason is to describe the regularizing part of the capillarity,
it will allow us to extract smooting effect and so to treat the non linear terms in distribution sense.
We recall that:
\begin{equation}
{\rm div}K
=\n\big(\rho\kappa(\rho)\D\rho+\frac{1}{2}(\kappa(\rho)+\rho\kappa^{'}(\rho))|\n\rho|^{2}\big)
-{\rm div}\big(\kappa(\rho)\n\rho\otimes\n\rho\big).
\label{divKa}
\end{equation}
and as for all $f\in C^{0}_{\infty}$:
$$\D f(\rho)=f^{'}(\rho)\D\rho+f^{''}(\rho)|\n\rho|^{2}.$$
We get then:
$${\rm div}K=\n(\D f(\rho)-\frac{1}{2}(\kappa(\rho)+\rho\kappa^{'}(\rho))|\n\rho|^{2})-{\rm div}(\kappa(\rho)
\n\rho\otimes\n\rho).$$
with $f^{'}(x)=x\kappa(x)$. It gives in particular the estimates (\ref{divK1}).
\subsection{Inequality energy estimates}
We are interested here in derivating bounds estimates on the system (\ref{3systeme}). We have to multiply momentum equation by $u$ and integrate over the time and the space.
We concentrate us only one term:
$$
\begin{aligned}
\int_{\R}\g {\rm div}K\cdot u \,dx\,dt=\int_{\R}\g&\biggl(\rho\kappa(\rho)\p_{i}\D\rho+\rho\kappa^{'}(\rho)\p_{i}\rho\D\rho+\frac{1}{2}\rho\kappa^{''}(\rho)
\p_{i}\rho|\n\rho|^{2}\\
&+\rho\kappa^{'}(\rho)\p_{j}\rho\p_{ij}\rho\biggl)u_{i}\,dx\,dt
\end{aligned}
$$
Next we have:
\begin{equation}
\begin{aligned}
&\int_{\R}\g\big(\frac{1}{2}\rho\kappa^{''}(\rho)
\p_{i}\rho|\n\rho|^{2}+\rho\kappa^{'}(\rho)\p_{j}\rho\p_{ij}\rho\big)u_{i}\,dx\,dt=\\
&-\hspace{1cm}\frac{1}{2}\int_{\R}\g \kappa^{'}(\rho)|\n\rho|^{2}u\cdot\n\rho\,dx\,dt-\frac{1}{2}\int_{\R}\g \rho\kappa^{'}(\rho)|\n\rho|^{2}{\rm div}u\,dx\,dt,
\end{aligned}
\label{appen1}
\end{equation}
and:
\begin{equation}
\begin{aligned}
&\int_{\R}\g\big(\rho\kappa(\rho)
\p_{i}\D\rho+\rho\kappa^{'}(\rho)\D\rho\p_{i}\rho\big)u_{i}\,dx\,dt=\\
&\hspace{1cm}-\int_{\R}\g \kappa(\rho)\D\rho\,u\cdot\n\rho\,dx\,dt-\int_{\R}\g \rho\kappa(\rho)\D\rho\,{\rm div}u\,dx\,dt.
\end{aligned}
\label{appen2}
\end{equation}
By mass equation we have:
\begin{equation}
 \rho{\rm div}u+u\cdot\n\rho=-\p_{t}\rho.
\label{appentrans}
\end{equation}
In using (\ref{appen1}), (\ref{appen2}) and (\ref{appentrans}) we get finally:
$$
\begin{aligned}
\int_{\R}\g {\rm div}K\cdot u \,dx\,dt=&\int_{\R}\g\kappa(\rho)\D\rho\p_{t}\rho\,dx\,dt+\frac{1}{2}
\int_{\R}\g\kappa^{'}(\rho)|\n\rho|^{2}\p_{t}\rho \,dx\,dt,\\
=&\frac{1}{2}\frac{\p}{\p t}\int_{\R}\g\kappa(\rho)|\n\rho|^{2} \,dx\,dt.
\end{aligned}
$$
The other terms are classical.


\begin{thebibliography}{}
\bibitem{3A}
D-. M. ANDERSON, G-.B McFADDEN and A- .A. WHELLER. Diffuse-interface
methods in fluid mech. \textit{In Annal review of fluid mechanics}, Vol. 30,
pages 139-165. Annual Reviews, Palo Alto,
CA, 1998.
\bibitem{3BCD}
H. BAHOURI , J.-Y. CHEMIN and R. DANCHIN. Fourier Analysis and Nonlinear Partial Differential Equations, \textit{Springer}, to appear.
\bibitem{3BDL}
D. BRESCH, B. DESJARDINS and C-. K. LIN, On some compressible fluid
models: Korteweg,lubrication and shallow water systems.  \textit{Comm.
Partial Differential Equations}, 28(3-4), 843-868, 2003.
\bibitem{3C}
J-. W. CAHN and J-. E. HILLIARD, Free energy of a nonuniform system, I.
Interfacial free energy,  \textit{J. Chem. Phys.} 28, (1998), 258-267.
\bibitem{3DD}
R. DANCHIN and B. DESJARDINS, Existence of solutions for
compressible fluid models of Korteweg type,  \textit{Annales de l'IHP, Analyse
non
lin\'eaire} 18, 97-133, (2001).
\bibitem{3DS}
J-. E. DUNN and J. SERRIN, On the thermomechanics of interstitial
working ,  \textit{Arch. Rational Mech. Anal.}, 88(2), (1985), 95-133.
\bibitem{3F}
E. Feireisl, Dynmamics of Viscous Compressible Fluids-Oxford Lecture
Series in Mathematics and its Applications-26 (2004).
\bibitem{3GP}
M-. E. GURTIN, D. POLIGONE and J. VINALS, Two-phases binary fluids and
immiscible fluids described by an order parameter,  \textit{Math. Models
Methods Appl. Sci.}, 6(6), (1996), 815-831.
\bibitem{3H}
B. HASPOT, Existence of strong solution for non isothermal Korteweg model,  \textit{to appear in Annales Blaise Pascal, 16, 2009}.
\bibitem{3H1}
H. HATTORI and D. LI, The existence of global solutions to a fluid
dynamic model
for materials for Korteweg type,  \textit{J. Partial Differential Equations}, 9(4), (1996), 323-342.
\bibitem{3H2}
H. HATTORI and D. LI, Global Solutions of a high-dimensional system
for
Korteweg materials,  \textit{J. Math. Anal. Appl.}, 198(1), (1996), 84-97.
\bibitem{3JL}
D. JAMET, O. LEBAIGUE, N. COUTRIS and J-. M. DELHAYE, The second
gradient method for the direct numerical simulation of liquid-vapor
flows
with phase change.  \textit{J. Comput. Phys.}, 169(2), 624--651, (2001).
\bibitem{3K}
D-. J. KORTEWEG. Sur la forme que prennent les \'equations du
mouvement des fluides si l'on tient compte des forces capillaires
par des variations de densit\'e.  \textit{Arch. N\'eer. Sci. Exactes S\'er.
II}, 6, 1-24, 1901.
\bibitem{3MK}
M. KOTSCHOTE.  Strong solutions for a compressible fluid model of Korteweg type.  \textit{Annales de l'Institut Henri Poincare (C) Non Linear Analysis},
Volume 25, Issue 4, July-August 2008, 679-696.
\bibitem{3L2}
P.-L. LIONS, Mathematical Topics in Fluid Mechanics, Vol 2,
Compressible models,  \textit{Oxford University Press}, (1996)
\bibitem{3MV}
A. MELLET and A. VASSEUR, On the isentropic compressible Navier-Stokes
equation, Arxiv preprint math.AP/0511210, 2005 - arxiv.org
\bibitem{3R}
J-. S. ROWLINSON, Translation of J.D van der Waals "The thermodynamic
theory of capillarity under the hypothesis of a continuous variation
of density". \textit{J. Statist. Phys.}, 20(2), 197-244, 1979.
\bibitem{3RS}
T. RUNST and W. SICKEL, Sobolev spaces of fractional order,
Nemytskij operators, and nonlinear partial differential
equations,\textit{ volume 3 of De Gruyter series in nonlinear analysis and applications.} Berlin 1996.
\bibitem{3Sim}
J. SIMON. Compact sets in the space $L^{p}(0,T;B)$.  \textit{Ann. Mat. Pura Appl.}, 146, 65-96, 1987.
\bibitem{3TN}
C. TRUEDELLAND and W. NOLL. The nonlinear field theories of mechanics.
 \textit{Springer-Verlag, Berlin, second edition}, 1992.
\end{thebibliography}
\end{document}